\documentclass[a4paper,11pt,english]{smfart}
 \usepackage{stmaryrd} 
\usepackage[applemac]{inputenc}
\usepackage[english]{babel}
\usepackage{amsfonts,mathrsfs}
\usepackage{amsmath}
\usepackage{amsthm}
\usepackage{hyperref} 
\usepackage{latexsym}
\usepackage{array}
\usepackage{amssymb}
\usepackage{enumerate}
\usepackage{smfthm}
\usepackage{graphicx}
\usepackage{color}
\newcommand{\ligne}{\vspace{1\baselineskip}}
\newcommand{\ph}{\phantomsection}
\newcommand{\cal}{\mathcal}
\newcommand{\di}{\displaystyle}

\newcommand{\R}{\mathbb  R}
\newcommand{\C}{\mathbb  C}
\newcommand{\N}{\mathbb  N}
\newcommand{\K}{\mathbb  K}

\renewcommand{\P}{\mathbb{P}}

\newcommand{\W}{  \mathcal{W}   }

\newcommand{\eps}{\varepsilon}
\renewcommand{\epsilon}{\varepsilon}
\newcommand{\e}{  \text{e}   }

\newcommand{\Z}{  \mathbb{Z}   }

\newcommand{\E}{\mathbb{E} }

\renewcommand{\H}{  \mathcal{H}   }

\newcommand{\dis}{\displaystyle}
\newcommand{\h}{  h }

\newcommand{\om}{  \omega   }
\newcommand{\ov}{  \overline  }
\renewcommand{\a}{  \alpha   }

\newcommand{\s}{  \sigma   }

\renewcommand{\phi}{  \varphi  }
\renewcommand{\L}{  \mathcal{L}   }

\newcommand{\<}{  \langle   }

\renewcommand{\>}{  \rangle   }

\renewcommand{\Im}{  \text{Im\,} }
\numberwithin{equation}{section}
\theoremstyle{plain}

\newtheorem*{acknowledgements}{Acknowledgements}

\def\beq{\begin{equation}}   \def\eeq{\end{equation}}
\def\bea{\begin{eqnarray}}  \def\eea{\end{eqnarray}}

\renewcommand{\theequation}{\thesection.\arabic{equation}}
\newcounter{hran} \renewcommand{\thehran}{\thesection.\arabic{hran}}

\def\bmini{\setcounter{hran}{\value{equation}}
    \refstepcounter{hran}\setcounter{equation}{0}
    \renewcommand{\theequation}{\thehran\alph{equation}}\begin{eqnarray}}

\def\bminiG#1{\setcounter{hran}{\value{equation}}
\refstepcounter{hran}\setcounter{equation}{-1}
\renewcommand{\theequation}{\thehran\alph{equation}}
\refstepcounter{equation}\label{#1}\begin{eqnarray}}

\author{ Didier Robert}
\address{Laboratoire de Math\'ematiques J. Leray, UMR  6629 du CNRS, Universit\'e de Nantes, 
2, rue de la Houssini\`ere,
44322 Nantes Cedex 03, France}
\email{didier.robert@univ-nantes.fr}
\author{ Laurent Thomann }
\address{Laboratoire de Math\'ematiques J. Leray, UMR  6629 du CNRS, Universit\'e de Nantes, 
2, rue de la Houssini\`ere,
44322 Nantes Cedex 03, France}
\email{laurent.thomann@univ-nantes.fr}

 \title[Random  weighted Sobolev inequalities]{Random  weighted Sobolev inequalities and application to  quantum  ergodicity}
 
\begin{document}
\frontmatter

 \begin{abstract}
 This paper is a continuation of \cite{PRT1} where we  studied  a randomisation method based on the Laplacian with harmonic potential. Here we extend  our previous results to the case of  any polynomial  and confining potential $V$ on $\R^d$.  We construct measures, under concentration  type assumptions, on the support of which we  prove optimal weighted Sobolev estimates on~$\R^d$. This construction relies on accurate estimates on the spectral function in a non-compact configuration space. Then we prove  random quantum ergodicity results  without specific assumption on the classical dynamics.  Finally, we prove that almost all basis of Hermite functions is quantum uniquely ergodic.
  \end{abstract}
\thanks{D. R.  was partly supported  by the  grant  ``NOSEVOL''   ANR-2011-BS01019 01.\\L.T. was partly supported  by the  grant  ``HANDDY'' ANR-10-JCJC 0109 \\
 \indent \quad and by the  grant  ``ANA\'E'' ANR-13-BS01-0010-03.}

\maketitle
\mainmatter

 \section{Introduction and results}
 \subsection{Introduction} 
 The aim of this paper is to prove stochastic properties of weighted Sobolev spaces associated to the 
 Schr\"odinger operator
 $$\hat H=\hat H_V = -\Delta + V,$$
  in  $L^2(\R^d)$ where $V$ is  a real confining  polynomial potential.
  We will see that most of the results  obtained  in~\cite{PRT1} for the harmonic oscillator ($V(x)=\vert x\vert^2$) 
   have a natural extension in this context. In the case of the harmonic oscillator we can take profit of many explicit formulas and both the spectrum and the eigenfunctions are well understood. This is not true anymore for more general potentials even if many asymptotic properties are known (see \cite{hero, iv, mo, ro1}). The results obtained in this paper combine accurate spectral  results for the operator $\hat H_V$  with probabilistic methods.    Our work is  inspired by similar properties proved before by Burq and Lebeau \cite{bule} (see also Shiffman and Zelditch~\cite{ShZe})
   concerning the Laplace operator on compact Riemaniann  manifolds.   When the degree of $V$ tends to infinity, we formally recover the estimates of  \cite{bule}. \ligne
   
   As a consequence of our probabilistic estimates, we obtain quantum ergodicity results for $\hat H$ under general probability laws satisfying the Gaussian concentration of measure property. These results extend previous ones \cite{zel1, bule} obtained on compact manifolds for the Gaussian measure.  
   Furthermore, our analysis allows to prove that almost all orthonormal bases of $L^{2}(\R^{d})$  of Hermite functions are unique quantum ergodic (see Theorem \ref{QUE} below).
   \\[2pt]

 \indent {\bf Spectral theory context :}
Let $d\geq 2$ and $k\geq 1$ an integer number and write $\<x\>=(1+\vert x\vert ^2)^{1/2}$. In the sequel, we will consider a potential $V$   satisfying\vspace{.2cm}
\vspace{.2cm}
    \begin{enumerate} 
 \item[(A1)]\qquad   The potential $V$  is  an elliptic  polynomial in $\R^d$, which means that  we have: $V= V_0 +V_1$ where $V_0$ is an homogeneous polynomial of degree $2k$ with ($V_0(x)>0$ if $x\neq 0$)
    and  $V_1(x) $ is a polynomial  of degree $\leq 2k-1$.
     \end{enumerate}
 \vspace{.2cm}
In particular for $\vert x\vert$ large we have $V(x) \approx \<x\>^{2k}$. Under these assumptions, $V$ is bounded from below;  by adding a constant to $V$, we can assume that~$V$ is nonnegative. Examples are $V(x) = \vert x\vert^{2k}$  for $k\in \N^{*}$  or $V(x) = \vert x\vert^{4} + \alpha\vert x\vert^2$,
$\alpha\in\R$.
\medskip

Notice that we could have worked with more general smooth and confining potentials. We have chosen to restrict ourselves to assumption (A1) in order to lighten the proofs and to get explicit exponents.
\medskip

Under these assumptions, $\hat H_V$ has a  compact resolvent  and its spectrum 
 consists in a non decreasing sequence  $\{\lambda_j\}_{j\geq 1}$ 
where each eigenvalue $\lambda_j$ is repeated according to its multiplicity and satisfies
 $\di{\lim_{j\rightarrow +\infty}\lambda_j = +\infty}$.  Let $\{\varphi_j\}_{j\geq 1}$ be an orthonormal basis in  $L^2(\R^d)$  of eigenfunctions of $\hat H_V$,
 $\hat H_V\varphi_j = \lambda_j\varphi_j$. Then one can show that  $\varphi_j$
is in the  Schwartz space  ${\cal S}(\R^d)$    thanks to Agmon estimates. 

For any  interval  $I$ of  $\R$  we denote by  $\Pi_H(I)$  the spectral projector of  $\hat H$ on  $I$.
 The range  ${\cal E}_H(I)$ of  $\Pi_H(I)$ is spanned  by  $\{\varphi_j; \lambda_j\in I\}$ and $\Pi_H(I)$ 
has an  integral kernel given by 
 \begin{equation}\label{def.pi}
 \pi_H(I;x,y) = \sum_{[j:\lambda_j\in I]}\varphi_j(x)\overline{\varphi_j(y)}.
 \end{equation}
If  $I = I_\lambda : =]0, \lambda]$, we denote by  $\Pi(\lambda)  = \Pi(I_\lambda)$, 
${\cal E}_H(\lambda)= {\cal E}_H(I_\lambda)$ and $\pi_H(\lambda; x,y)=\pi_H(I; x,y)$.\\
We have a "soft"  Sobolev inequality (see e.g. \cite[Lemma 3.1]{PRT1} for a proof):
\beq\label{sobg}
 \vert u(x)\vert \leq (\pi_H(I; x,x))^{1/2}\Vert u\Vert_{L^2(\R^d)},\;\; \forall u\in {\cal E}_H(I).
\eeq
 In $\R^d$ the spectral function   $\pi_{H_V}(\lambda; x,x)$ is fast decreasing for  $\vert x\vert\rightarrow +\infty$  so
  it is natural to replace the~$L^\infty$ norm   $\di{\Vert u\Vert_\infty = \sup_{x\in\R^d}\vert u(x)\vert}$
  by  the weighted norms 
  $\di{\Vert u\Vert_{\infty, s} = \sup_{x\in\R^d}\<x\>^s\vert u(x)\vert}$.\ligne

The spectral theory for  $\hat H$ has been studied in    \cite{hero, mo, iv}. In particular a Weyl asymptotics was proven:
     \beq\label{Weylhom}
     N_H(\lambda) =   (2\pi)^{-d} \int_{[\vert\xi\vert^2 + V(x)\leq \lambda]}d\xi dx + \mathcal{O}(\lambda^{\frac{(d-1)(k+1)}{2k}}),\; \lambda \rightarrow +\infty,
     \eeq
and using the assumptions (A1) on $V$ we get   $ N_H(\lambda) \approx \lambda^{\frac{d(k+1)}{2k}}$. 
We then consider the normalized Hamiltonian  
 \begin{equation}  \label{defHnor}
 \hat H_{nor}: =\hat H^{\frac{k+1}{2k}},
  \end{equation}
  so  that we have 
 \begin{equation*}     
      N_{H_{nor}}(\lambda) \approx \lambda^d.
      \end{equation*}
The Sobolev spaces associated with $\hat H$
      are here defined as follows. Let  $s\geq 0$, $p\in[1, +\infty]$,
      \beq\label{sobsp}
      {\cal W}_k^{s,p}:={\cal W}_k^{s,p}(\R^d):=\big\{u\in L^p(\R^d),\; \hat H_{nor}^su\in L^p(\R^d)\big\},\; \Vert u\Vert_{s,p} = \Vert\hat H_{nor}^su\Vert_{L^p(\R^d)}.
       \eeq
       
       More explicitly, for $1<p<+\infty$ we can prove  (see \cite[Lemma~2.4]{YajimaZhang2}) 
        \begin{equation*}
          {\cal W}_k^{s,p}:={\cal W}_k^{s,p}(\R^d):=\big\{u\in L^p(\R^d),\; (-\Delta)^{s\frac{k+1}{2k}}u\in L^p(\R^d)\;\; \mbox{and} \;\;|x|^{s(k+1)}u\in L^p(\R^d)\big\}.
              \end{equation*}

       The Hilbert Sobolev spaces are denoted by  ${\cal H}_k^s =   {\cal W}_k^{s,2}$.  
        The ${\cal H}_{k}^s$-norm is equivalent to the spectral norm
  \begin{equation*}   
        \Vert u\Vert_{{\cal H}_{k}^s} = \Big(\sum_{j\in\N}\lambda_j^{\frac{(k+1)s}{k}}\vert c_j\vert^2\Big)^{1/2},\;\; {\rm when}\;\;\;\;\; u(x)=\sum_{j\in\N}c_j\varphi_j(x).
       \end{equation*}
 \begin{rema}
 We stress that with the definition \eqref{sobsp}, $\hat H_{nor}$ is considered  as an order 1 operator. This convention is different from the one used in \cite{PRT1}, where the definitions of the Sobolev spaces were based on $\hat{H}$ instead of $\hat H_{nor}$.
 \end{rema}

 Let $d\geq 2$.   For $h>0$, we define the spectral  interval   $I_h = [\frac{a_h}{h}, \frac{b_h}{h}[$  and we assume that $a_{h}$ and $b_{h}$ satisfy, for some $a,b,D>0, \delta\in[0, 1]$, 
\begin{equation}\label{spect1}
 \lim_{h\rightarrow 0}a_h =a,  \;\;\; \lim_{h\rightarrow 0}b_h = b,\quad 0<a \leq b\quad \text{and}\quad  b_h -a_h \geq Dh^\delta.
\end{equation}
with any  $D>0$ if $\delta<1$ and $D$ some large constant if $\delta=1$.
We denote by $N_h$ the number with multiplicities of eigenvalues of $\hat{H}_{nor}$ in $I_h$.  By \eqref{Weylhom} we have
\begin{equation}\label{eqN}
N_{h}\sim (2\pi)^{-d} \int_{\frac{a_{h}}{h}\leq \big(\vert\xi\vert^2 + V(x)\big)^{{(k+1)}/{(2k)}}<\frac{b_{h}}{h}}d\xi dx.
\end{equation}
Then, under the previous assumptions, $\dis \lim_{h\rightarrow 0}N_h = +\infty$. 

 In the sequel, we write $\Lambda_h = \{j\geq 1, \;\lambda^{\frac{k+1}{2k}}_j\in I_h\}$ and ${\cal E}_h = {\rm span}\{\phi_{j},\;j\in \Lambda_h\}$, so that  $N_h  = \#\Lambda_h={\rm dim}\,{\cal E}_h$. We denote by ${\bf S}_h=\big\{u\in \cal{E}_h\,:\;\; \Vert u\Vert_{L^2(\R^d)} =1\big\}$ the unit sphere  of~$\cal{E}_h$.\ligne

We will consider sequences of complex numbers  $\gamma = (\gamma_{n})_{n\in \N}$  so that  there exists $K_{0} >0$ 
 \beq\label{condi1}
 \vert\gamma_n\vert^2 \leq \frac{K_0}{N_h}\sum_{j\in\Lambda_h}\vert\gamma_j\vert^2, \quad \forall n\in\Lambda_h,\;\;\forall h\in ]0, 1].
\eeq
This condition  means that on each level of energy $\lambda_n$, $n\in\Lambda_h$, one  coefficient $|\gamma_{k}|$ cannot be much larger  than the others. Sometimes, in order to prove lower bound estimates, we will need the stronger condition ($K_{1}>0$)
\beq\label{condi0}
 \frac{K_1}{N_h}\sum_{j\in\Lambda_h}\vert\gamma_j\vert^2 \leq \vert\gamma_n\vert^2 \leq \frac{K_0}{N_h}\sum_{j\in\Lambda_h}\vert\gamma_j\vert^2, \quad \forall n\in\Lambda_h,\;\;\forall h\in ]0, 1].
\eeq
  
 This condition which was suggested to us by Nicolas Burq,  means that on each level of energy $\lambda_n$, $n\in\Lambda_h$,  the coefficients~$|\gamma_{k}|$ have almost the same size. For instance \eqref{condi0} holds if there exists $(d_{h})_{h\in]0,1]}$ so that $\gamma_{n}=d_{h}$ for all $n\in \Lambda_h$. These sequences will be used to built  random wave packets
 $u_\gamma = \sum \gamma_jX_j\varphi_j$ for random variables $\{X_j\}_{j\geq 1}$.
 \ligne

 {\bf Stochastic context :}
Consider a   probability space
$(\Omega, {\cal F}, {\P})$  and  let $\{X_n,\;n\geq1\}$ an i.i.d  sequence of real or complex random variables with joint  distribution   $\nu$. In the sequel, we assume that  the r.v.\,$X_{n}$ is centred ($\E(X_{n})=0$) and normalised 
($\E( \vert X_{n}\vert^{2})=1$).
 For our results we will assume that $\nu$ satisfies   a Gaussian concentration condition  in the sense of the following definition:

 \begin{defi}\ph 
We say that a probability measure   $\nu$ on $\K$  ($\K=\R,  \C$)  satisfies the concentration of measure property if there exist  constants $c,C>0$ independent of $N\in \N$ such that for all Lipschitz and convex function $F : \K^{N}\longrightarrow \R$
\begin{equation}\label{Blip}
\nu^{\otimes N}\Big[\,X\in \K^{N}\;:\;\big|F(X)-\E(F(X))\big|\geq r \,\Big]\leq c\,e^{-\frac{Cr^{2}}{ \|F\|^{2}_{Lip}}},\quad\forall r>0, 
\end{equation}
where $\|F\|_{Lip}$ is the best constant so that $\dis |F(X)-F(Y)|\leq \|F\|_{Lip}\|X-Y\|_{\ell^{2}}$.
\end{defi}
For instance,  $\nu$  can be the   standard normal law,  a Bernoulli law or any law with bounded support (see \cite[Section 2]{PRT1}).  For details on this notion, see Ledoux \cite{led}.

Observe that if \eqref{Blip} is satisfied, that for $\eps>0$ small enough, for $j\geq 1$ (see \cite[(2.3)]{PRT1})
\begin{equation}\label{quad}
\int_{\K}\e^{\eps X^{2}_{j}}\text{d}\nu(x)<+\infty.
\end{equation}
~

If $(\gamma_{n})_{n\in \N}$ satisfies \eqref{condi0} and $\nu$ satisfies  \eqref{Blip}, we  define the random vector  in ${\cal E}_h$ 
\begin{equation*} 
v_{\gamma}(\omega) : = v_{\gamma,h} (\omega)= \sum_{j\in\Lambda_{h}}\gamma_jX_j(\omega)\phi_{j},
 \end{equation*}
 and we  define a probability  ${\bf P}_{\gamma,h}$ on ${\bf S}_h$ by: for all measurable and bounded function $f:{\bf S}_h\longrightarrow \R$
 \begin{equation}\label{ProbaS}
  \int_{{\bf S}_h}f(u)\text{d}{\bf P}_{\gamma,h}(u) =  \int_\Omega f\left(\frac{v_\gamma(\omega)}{\Vert v_\gamma(\omega)\Vert_{L^{2}(\R^{d})}}\right)\text{d}\P(\omega).
 \end{equation}
In other words,  ${\bf P}_{\gamma,h}$ is the push forward probability  measure obtained  by projection on the unit sphere of ${\cal E}_h$
 of the law of  the random vector $v_\gamma$.  In the  isotropic case ($\gamma_j = \frac{1}{\sqrt N}$ for all $j\in\Lambda_{h}$),  we denote this probability by ${\bf P}_{h}$. We can check (see e.g. \cite[Appendix C]{PRT1}) that in the isotropic case and in the particular case where  $\nu=\mathcal{N}_{\C}(0,1)$ or $\nu=\mathcal{N}_{\R}(0,1)$, then ${\bf P}_{h}$ (which we will denote by ${\bf P}^{(u)}_{h}$)   is the uniform probability on~${\bf S}_h$. 
 
 Notice that if the probability law $\nu$ is rotation invariant on $\C$ then the  probabilities ${\bf P}_{\gamma,h}$ are invariant by the 
 Schr\"odinger  linear flow ${\rm e}^{-ith^{-1}\hat H}$.
 
 \subsection{Main results of the paper}
Let us state now the main results proved in this paper.

  \subsubsection{\bf Estimates  for frequency localised functions} 
Our first result shows that for the elements in  the support of ${\bf P}_{\gamma,h}$, the Sobolev estimates are  better than usual. This  phenomenon   was studied  in \cite{bule}  on compact Riemannian manifolds  and in \cite{PRT1} for the harmonic oscillator. Recall the definition \eqref{sobsp} of the Sobolev spaces $\W_{k}^{s,p}(\R^{d})$, then 
  \begin{theo} \ph\label{thm1}
Let $d\geq 2$ and assume that $0\leq \delta<2/3$ in \eqref{spect1}.  Consider a potential $V$ which satisfies conditions $(A1)$. Suppose that $\nu$ satisfies the concentration of measure property \eqref{Blip}. Then there exist $0<C_0 < C_1$,    $c_1>0$ and $h_0>0$ such that for all $h\in]0, h_0]$.
\begin{equation} \label{prest1}
{\bf P}_{\gamma,h}\left[u\in{\bf S}_h : C_0\vert\log h\vert^{1/2}  \leq 
\Vert u\Vert_{\W_k^{\frac{d}{2(k+1)},\infty}(\R^{d})}  
\leq C_1\vert\log h\vert^{1/2}\,  \right] \geq 1- h^{c_1}.
\end{equation}
Moreover the estimate  from above holds true  for  any $0\leq \delta \leq 1$ ($D$ is large enough for $\delta=1$)
and $\gamma$ satifying (\ref{condi1}).
\end{theo}
Recall that on the support of ${\bf P}_{\gamma,h}$   there exist $0<C_2 < C_3$, so that for all $u\in {\bf S}_h$, and $s\geq 0$
\begin{equation}\label{equi}
C_2h^{-s}  \leq \Vert u\Vert_{\H_{k}^{s}(\R^{d})}\leq C_3 h^{-s},
\end{equation}
which shows that there is no random smoothing in the $L^2$ Sobolev scale spaces (see also Theorem\,\ref{nors}).
But, if $\hat H=-\Delta+V$ is considered as an operator of order 2, then using \eqref{defHnor}, Theorem \ref{thm1} shows a gain of  $d/(2k)+d/2$ derivatives in $L^{\infty}(\R^{d})$
comparing with the deterministic Sobolev embedding.  We shall see  that this is really a stochastic smoothing property  which has nothing to do with the usual Sobolev embedding.  

For $k=1$ we recover \cite[Theorem 1.1]{PRT1} in the case of the harmonic oscillator. In this particular case, we moreover construct Hilbertian bases of $L^{2}(\R^{d})$ of eigenfunctions of $\hat H$ which enjoy good decay properties (see \cite[Theorem 1.3]{PRT1}). For a general potential we can extend this result, with  analogous bounds, but for quasimodes instead of exact eigenfunctions. 
 
  When $k\to +\infty$, the operator $\hat H$ seems to  behave like the Laplacian on a compact manifold. In this case we get  estimates close to  the estimates of \cite[Théorème 5]{bule}.
  
The proof of Theorem \ref{thm1} relies on a good understanding of the spectral function associated with~$\hat H$, in particular on two weighted estimates in $L^{p}$: The first is a semi-classical uniform upper bound obtained by Koch-Tataru-Zworski \cite{KTZ}, while the second estimate is a lower bound which is proved in the Appendix. These results are combined with probabilistic techniques as in \cite{bule,PRT1}.\ligne

  \subsubsection{\bf Random quantum ergodicity} 
Our second result concerns semi-classical measures  associated with families of states 
$h\mapsto u_h$, $h\in]0, 1]$, $\Vert u_h\Vert_{L^2(\R^d)}=1$. This development is inspired from \cite{bule}, and can be interpreted as a kind of random-quantum-ergodicity (see \cite{zel1}).

Assume that  $J_h:=hI_h= ]a_h, b_h]$ is such that 
    \begin{equation*}
        \lim_{h\rightarrow 0}a_h=  \lim_{h\rightarrow 0}b_h=\eta>0\qquad  \text{and}\qquad 
    \lim_{h\rightarrow 0}\frac{b_h -a_h}{h}=+\infty. 
    \end{equation*}
 
       Denote by $L_\eta$ the Liouville measure  associated with the classical Hamiltonian 
    $H_0(x,\xi) = \frac{\vert \xi\vert^2}{2} + V_0(x)$ (recall that $V_{0}$ is defined in Assumption $(A1)$). Then $ L_\eta$ is given by 
    $$
    L_\eta(A) = C_\eta\int_{[H_0(z)=\eta]}\frac{A(z)}{\vert\nabla H_0(z)\vert}d\Sigma_\eta(z)
    $$
    where $d\Sigma_\eta$ is the Euclidean measure on the hyper surface $\Sigma_\eta:=H_0^{-1}(\eta)$ and  $C_\eta>0$ is a normalization constant such that $L_\eta$ is a probability measure on $\Sigma_\eta$.
    
We denote by $S(1,k)$  the class of symbols   such that $A\in C^\infty(\R^{2d})$ and $A$ is quasi-homogeneous of degree 0 outside a small neighborhood of $(0,0)$ in $\R_x^d\times\R_\xi^d$:
 $A(\lambda x, \lambda^k\xi) = A(x,\xi)$ for every $\lambda\geq 1$ and $\vert(x,\xi)\vert \geq \eps$.
 
For $A\in S(1,k)$ let us denote by $\hat A$ the Weyl quantization of $A$ (here $h=1$).
Notice that if $\Vert u_h\Vert_{L^2(\R^d)} =1$ then $A\mapsto \<u_h, \hat Au_h\>$ defines a semiclassical measure on $\Sigma_\eta$ (see e.g. \cite{bbk}).  
\ligne 

 Then we have
 \begin{theo}\ph\label{rqer}
       Consider a potential $V$ which satisfies conditions (A1). Assume that we are in the isotropic case  ($\gamma_j = \frac{1}{\sqrt N_h}$ for all $j\in\Lambda_{h}$), and that $\nu$ satisfies  the concentration of measure property~\eqref{Blip}.
       Then there exist $c,C>0$ so that for all $r\geq 1$ and $A\in S(1,k)$ with $\dis {\max_{(x,\xi)\in \Sigma_{1}}|A(x,\xi)|\leq 1}$
  \begin{equation}\label{GD}
  {\bf P}_{h}\Big[ \,u\in {\bf S}_{h} : |\<u, \hat A u\>-L_{\eta}(A)|>r   \Big]\leq C\e^{-cN_hr^{2}}.
  \end{equation}    
       \end{theo}
       By homogeneity we can also get an estimate for any $A\in S(1,k)$.
       
        We are inspired from    Burq-Lebeau \cite[Théorème 3]{bule} who  obtained a similar result for the Laplacian on a compact manifold. A 
modification of their proof allows to consider more general random variables satisfying the Gaussian concentration
assumption instead of the uniform law.

   The  result of Theorem \ref{rqer} can be related to  quantum ergodicity (see \cite{zeq, cdv, hmr}) which concerns the semi-classical
     behavior of  $\<\varphi_j, \hat A\varphi_j\>$ when the classical flow is ergodic on the energy hyper surface\,$\Sigma_\eta$. Then, for "almost all" eigenfunctions $\varphi_j$,
     we have
      $\<\varphi_j, \hat A\varphi_j\>\stackrel{j\rightarrow +\infty}{\longrightarrow} L_\eta(A)$.
      The meaning of Theorem\,\ref{rqer} is that we have 
      $\<u, \hat Au\>\stackrel{h\rightarrow 0}{\longrightarrow} L_\eta(A)$
      for almost all $u\in{\bf S}_h$ such that all modes $(\varphi_j)_{j\in\Lambda_h}$ are   ``almost-equi-present" in $u$.  For related results see Zelditch \cite{zel2}.
      
       As a consequence of Theorem \ref{rqer}
      we easily get  that almost all bases  of Hermite functions are Quantum Uniquely Ergodic (see
      Theorem \ref{QUE} for a precise statement).

     From \eqref{GD} we directly deduce that there exists $C>0$ so that for all $p\geq 2$ and $h\in]0, 1]$,
    \begin{equation*}
    \big\|\<u, \hat A u\>-L_{\eta}(A)\big\|_{L^{p}_{{\bf P}_{h}}}\leq CN_h^{-1/2}\sqrt{p}.
    \end{equation*} 
Therefore, if one denotes by  
      $$
      u_h^\omega = \frac{1}{\sqrt{N_h}}\sum_{j\in\Lambda_h}X_j(\omega)\varphi_j,
      $$
      then we have  
      $\<u^\omega_h, \hat Au^\omega_h\>\stackrel{h\rightarrow 0}{\longrightarrow} L_\eta(A)$
     in $L^p(\Omega)$-norm with a remainder estimate.  For example, we can apply this result to    $\di{ u_h^\omega = \frac{1}{\sqrt{N_h}}\sum_{j\in\Lambda_h}(-1)^{\epsilon_j(\omega)}\varphi_j}$
      where $\{\epsilon_j\}_{j\geq 1}$ are i.i.d Bernoulli variables.
      
      \begin{rema}
      It is likely that our main results can be extended in the case where $V$ is a smooth (instead of polynomial) confining potential, at the cost of some technicalities. 
      \end{rema}
      
        \subsubsection{\bf Quantum unique ergodicity of bases of Hermite functions} 
          As another application of the ergodicity results (see Proposition \ref{scquerg}) we prove that  random orthonormal basis of eigenfunctions of the Harmonic oscillator $\hat H=-\Delta+|x|^{2}$ is  Quantum Uniquely Ergodic (QUE, according the terminology used in~\cite{zel2}
     and introduced in \cite{rusa}).
     
     Firstly, we assume that for all $j\in \Lambda_{h}$, $\gamma_{j}=N_h^{-1/2}$ and that $X_{j}\sim \mathcal{N}_{\C}(0,1)$, so that ${\bf P}_{h}:={\bf P}_{\gamma,h}$ is the uniform probability on ${\bf S}_{h}$. Set $\eta=2$, then the Liouville measure $L_{2}$ is the uniform measure on $\big\{(x,\xi)\in \R^{2d}:\; |x|^{2}+|\xi|^{2}=2\big\}=\sqrt{2}\,\mathbb{S}^{2d-1}$. We set $h_{j}=1/j$ with $j\in \N^{*}$, and  
 \begin{equation*}
 a_{h_{j}}=2+dh_{j},\quad b_{h_{j}}=2+(2+d)h_{j}.
 \end{equation*}
  Then   \eqref{spect1}  is satisfied with  $\delta=1$ and  $D=2$. In particular, each interval
   $$I_{h_{j}}=\Big[\,\frac{a_{h_{j}}}{h_{j}},\frac{b_{h_{j}}}{h_{j}}\Big[=[2j+d,2j+d+2[$$
    only contains  the  eigenvalue $2j+d$ with multiplicity $N_{h_{j}}\sim cj^{d-1}$, and  $\mathcal{E}_{h_{j}}$ is the corresponding eigenspace of the harmonic oscillator $H$. We can identify the space of the orthonormal basis of $\mathcal{E}_{h_{j}}$ with the unitary group $U(N_{h_{j}})$ and  we endow $U(N_{h_{j}})$ with its Haar probability measure $\rho_{j}$.  Then the space $\mathcal{B}$ of the  Hilbertian bases of  eigenfunctions of $H$ in $L^{2}(\R^{d})$ can be identified with 
 \begin{equation*}
 \mathcal{B}=\times_{j\in \N} U(N_{h_{j}}),
 \end{equation*}
 which can be endowed with the measure 
 \begin{equation*}
\dis \text{d} \rho=\otimes_{j\in \N}\,\text{d} \rho_{j}.
 \end{equation*}
Denote by $B=(\phi_{j,\ell})_{j\in \N, \,\ell\in  \llbracket 1,N_{h_j} \rrbracket } \in \mathcal{B}$ a typical orthonormal basis of $L^{2}(\R^{d})$ so that for all $j\in \N$, $(\phi_{j,\ell})_{\ell\in  \llbracket 1,N_{h_j} \rrbracket } \in U(N_{h_{j}})$ is an orthonormal basis of $\mathcal{E}_{h_{j}}$. 

\begin{theo} \ph \label{QUE} Let $\hat H$ be the harmonic oscillator. 
For $B\in\mathcal{B}$ and $A\in S(1,1)$  denote by
$$
D_j(B) = \max_{1\leq \ell\leq N_{h_j}}\big\vert\<\varphi_{j,\ell}, \hat A\varphi_{j,\ell}\> -L_{2}(A)\big\vert.
$$
Then we have
$$ \lim_{j\rightarrow +\infty}D_j(B) = 0,\;\; \rho-{\rm a.s\; \,on}\;\;\mathcal{B}.$$
In other words, a random orthonormal basis of Hermite functions is QUE.
  \end{theo}

In \cite{zel1} the author proved that on the standard sphere
      a  random orthonormal basis of eigenfunctions of the Laplace operator is ergodic.  In \cite{map} the author get  a unique quantum ergodicity  result   for compact  Riemannian manifolds.
      
       Using  estimates proved in \cite{bule}, an analogous result  to Theorem \ref{QUE} can be obtained for the Laplace operator  on Riemannian compact manifolds  with the same method. This holds true in particular for the sphere in any dimension $d\geq 2$  and more generally for   Zoll manifolds (in this last setting   a random orthonormal  basis  of quasi-modes is  obtained).

      \subsection{Plan of the paper} The rest of the paper is organised as follows.    In  Section \ref{Sect2} we state crucial estimates on the spectral function which are used in    Section \ref{Sect3}  for the proof of Theorem \ref{thm1}.  Section \ref{Sect4} is devoted to the proof of   Theorem \ref{rqer}. Finally, in the appendix, we prove some point-wise estimates for the spectral function.
      
      \begin{acknowledgements}
The authors benefited from discussions with  Nicolas Burq and  Aur\'elien Poiret.
\end{acknowledgements}

 
  \section{Estimates for the spectral function}\label{Sect2}
  
     Our goal here is to obtain  for the spectral function of operators $\hat H$ satisfying (A1),  analogous estimates as those proved in \cite{kar} and in \cite{KTZ}  for the harmonic oscillator. These estimates are more or less known, in particular  some are contained in \cite{kar}, but they  are  crucial for  proving our first main result,  so we state them in a form suitable for us. The difficulties to prove them increase with the parameter $\delta\in[0, 1]$   measuring the width of  the spectral window $I_h$.  Concerning upper bounds, the case $\delta=0$ is easy; $\delta <2/3$ can be reached using accurate properties of the semi-classical Weyl calculus; the case $\delta =1$ was  solved in \cite{KTZ}  using a semi-classical Strichartz estimate. However, for the lower bounds, we need to assume  $\delta < 2/3$.
     
     \subsection{An upper bound of the spectral function}       Here we cannot use explicit formulas like the Mehler formula for the harmonic oscillator. But it is possible to use a global pseudo-differential calculus to construct
  a parametrix for        $(\hat H -\lambda)^{-1}$ where   $\lambda\in\C$, $\lambda\notin[0, \infty[$  
       and to use some estimates proved in~\cite{ro1}. We  give in Appendix \ref{appA2} the 
       main steps to prove the following result.
              \begin{prop}\ph\label{esthom}
      Let us denote by  $K_H(t; x,y)$ the heat kernel of  $\hat H=-\Delta+ V$. Then  for every $N\geq 1$ there exists $C_N>0$  such that 
  \begin{equation*}
    K_H(t;x,x) \leq C_N\Big(t^{-d/2}\exp(-tV(x))  + 
       (1+\vert x\vert^{2k})^{-N}\Big) 
  \end{equation*}
          for every $t\in]0, 1]$ and  $x\in\R^d$.
    \end{prop}
   The bound stated in Proposition \ref{esthom} is relevant for  $t\to 0$. For better estimates in the regime $t>0$ we refer e.g. to \cite[Proposition 1]{Sikora} and references therein.\medskip
   
   We get easily the following estimate for the spectral function of  $\hat H$ defined in \eqref{def.pi}.
     \begin{coro} 
     With the notations of the previous Proposition we have
\begin{equation*}
      \pi_H(\lambda;x,x) \leq C_N\left(\lambda^{d/2} 
      \exp\big(-V(x)/{\lambda}\big)  + 
       (1+\vert x\vert^{2k})^{-N}\right),\;\;\forall \lambda\geq 1.
\end{equation*}
     In particular for every $\theta \geq 0$ we have 
\begin{equation*}
     \pi_H(\lambda;x,x) \leq C\lambda^{(d +\theta)/2} \<x\>^{-k\theta}
\end{equation*}
and for every  $u\in{\cal E}_H(\lambda)$, thanks to \eqref{sobg} we have 
\begin{equation}\label{sobhk}
   \<x\>^{k\theta/2}\vert u(x)\vert \leq C_\theta\lambda^{\frac{d +\theta}{4}}\Vert u\Vert_{L^2(\R^d)}.
\end{equation}
   \end{coro}  
      \begin{rema} 
      Recall  that $N_H(\lambda)$ (number of eigenvalues  $\leq \lambda$)
      is  of order $\lambda^{d_k}$  with  $d_k = \frac{d(k+1)}{2k}$. \\
     We shall see that  $\theta = \frac{d}{k}$   is specific  because we get from \eqref{sobhk}
\begin{equation*}
     \<x\>^{d/2}\vert u(x)\vert \leq C\lambda^{\frac{d_k}{2}}\Vert u\Vert_{L^2(\R^d)}.
\end{equation*}
  Notice that we also have 
       \begin{equation*}
              \pi_H(\lambda;x,x) \leq C\lambda^{d/2}.
       \end{equation*}
       \end{rema}
     Now using the definition of $h=\lambda^{-\frac{k+1}{2k}}$ we get that   
\begin{equation*}
 \<x\>^{k\theta/2}h^{\frac{k(d+\theta)}{2(k+1)}}\vert u(x)\vert \leq C\Vert u\Vert_{L^2(\R^d)}, \quad
\forall u\in {\cal E}_{H_{nor}}(h^{-1}).
\end{equation*}

\subsection{A uniform estimate for the spectral function}
  Now we shall state more refined estimates.  Let $V$ be a potential satisfying (A1),  then we have
  
    \begin{prop}\ph\label{kaest}
 For every $\delta\in[0, 1]$  and $C_0>0$,  there exists $C>0$ such that for every $\theta\in[0, \frac{d}{k}]$  and $1\leq r\leq \infty$  
 \beq\label{splp}
 \Vert\pi_{H_{nor}}(\lambda+\mu; x,x) - \pi_{H_{nor}}(\lambda; x,x)\Vert_{L^{r, k(r-1)\theta}(\R^d)} \leq C\lambda^\alpha
 \eeq
  for $\vert\mu\vert \leq C_0\lambda^{1-\delta}$, $\lambda \geq 1$ and  with  $\alpha = \frac{d}{k+1}(k + \frac{1}{r})-\delta +\frac{k\theta}{k+1}(1-\frac{1}{r})$,
    where $\Vert u\Vert_{L^{r,s}(\R^d)}^r = \int_{\R^d}\<x\>^s\vert u(x)\vert^rdx$.
  \end{prop}

  The proof Proposition~\ref{kaest} is a consequence of the following more general semiclassical result, which is proved in Appendix~\ref{A.3}. See also Ivrii \cite{iv} for related work and \cite[Lemma 3.5]{PRT1} for more details. For $0<\eps<1$, denote by ${V_{\eps}(x)=\eps^{2k}V(x/\eps)}$, where $V$ satisfies (A1) and consider the operator $\hat H_{\eps}(h)=-h^{2}\Delta+V_{\eps}$. Then  
  \begin{prop}\ph\label{scsfest} For every $\delta\in[0, 1]$,  every  $C_0 >0$
 there exist $\varepsilon_0>0$  and $C_1>0$  such that
   for every  $\nu\in]\nu_0-\varepsilon_0, \nu_0+\varepsilon_0[$,      $\vert\mu\vert \leq C_0h^{\delta-1}$,  $x\in\R^d$,  
   $h\in ]0, 1]$  and $\eps\in]0, 1]$ we have 
\begin{equation*} 
  \vert\pi_{H_\eps(h)}(\nu+\mu h;x,x) - \pi_{H_\eps(h)}(\nu;x,x)\vert \leq C_1h^{\delta -d}.
  \end{equation*}
   \end{prop}
   \begin{proof} See Appendix \ref{A.3}.
    \end{proof}
    
    \begin{rema}
    Actually, this result holds true for a larger class of potentials  $(V_{\eps})_{\eps\in(0,1]} \subset \mathcal{C}^{\infty}(\R^{d})$, which satisfy the following conditions:   Let $\nu_{0}\in \R$, $\eps_{0}>0$ and denote by $I_{0}=[\nu_{0}-\eps_{0},\nu_{0}+\eps_{0}]$. We assume that there exists a compact $K\subset \R^{d}$ so that,   uniformly in $\eps\in(0,1]$ 
    \begin{equation*}
      \left\{
      \begin{aligned}
         &\bullet \;\; V^{-1}_{\eps}(I_{0}) \subset K \\
 &\bullet \;\; |\partial^{\a}_{x} V_{\eps}(x)|\leq C_{\a} \;\;\text{for all}\;\; x\in K \;\;\text{and}\;\; \a\in \N^{d}\\
&\bullet \;\; \text{There exists}\;\;  \delta_{0}>0  \;\;\text{so that}\;\;  |\nabla V_{\eps}(x)|\geq \delta_{0}  \;\;\text{for all}\;\;   x\in K.
\end{aligned}
\right.
    \end{equation*}
   \end{rema}
   \begin{rema}
   As we can see in Appendix \ref{A.3}, the proof uses standard semiclassical arguments, with the semiclassical parameter $h\approx \lambda^{-\frac{k+1}{2k}}$,  where $\lambda$ is the natural energy parameter of our initial Hamiltonian $H_{V}$ (here $2k$ is the degree of homogeneity of $V$).
    
   	As already said, for $0\leq \delta < 2/3$ the Proposition can be proved using semiclassical pseudo-differential Weyl calculus.
	The general case $0\leq \delta \leq 1$ is more difficult\footnote{ We thank M. Zworski for having explained to us
	that the results of \cite{KTZ} give the $L^\infty$-estimate  for the spectral function with $\delta=1$,  in particular
at the turning points.}. We shall see in Appendix that it is a consequence of results proved in
	\cite{KTZ} extended to Hamiltonians depending smoothly on a parameter~$\eps$.
     \end{rema}

 \begin{proof}[Proof of Proposition \ref{kaest}]
We reduce  the problem  to a semi-classical estimate.\\
 Let us consider the eigenvalue problem
 \beq\label{vp1}
 (-\Delta_x +V(x))\varphi_j(x) = \lambda_j\varphi_j(x).
\eeq
Consider the $(\psi^{h}_{j})$ defined by 
\begin{equation}\label{defpsi}
\psi^h_j(y) = h^{-\frac{d}{2(k+1)}}\varphi_j(h^{-\frac{1}{k+1}}y).
\end{equation}
 Performing  the change of variable $x=yh^{-\frac{1}{k+1}}$, equation (\ref{vp1}) is transformed into
\begin{equation*} 
( -h^2\Delta_y + V_{\varepsilon}(y))\psi^h_j(y) = \lambda_j(h)\psi^h_j(y),
 \end{equation*}
 $\lambda_j(h) = \lambda_jh^{\frac{2k}{k+1}}$, and where $V_\varepsilon$ is the potential
 $$V_{\eps}(y)=\eps^{2k}V(\frac{y}\eps),$$
which is   smooth  in the parameter  $\varepsilon\in[0, 1]$ and is computed for $\varepsilon = h^{\frac{1}{k+1}}$.\\
So  $\lambda \leq \lambda_j^{\frac{k+1}{2k}} \leq
\lambda +\mu$
if and only if $\lambda h \leq \big(\lambda_j(h)\big)^{\frac{k+1}{2k}}\leq
\lambda h +\mu h$, therefore we can fix a large real parameter~$\tau$ and work with the semiclassical parameter 
$h = \frac{\tau}{\lambda}$ which is small for $\lambda$ large.   Denote by $\hat H_\eps(h) = -\h^2\Delta + V_\epsilon(x)$. 
   Then the potential
   $V_\epsilon$  satisfies the condition (A1) and we can apply Proposition~\ref{scsfest}. We deduce that for every $c>0$  and for every $\delta \leq  1 $ there exists $C>0$ such that
\begin{equation*} 
 \vert\pi_{H_{nor}}(\lambda+\mu;x,x) -  \pi_{H_{nor}}(\lambda;x,x)\vert \leq C\lambda^{d\frac{k}{k+1} -\delta}
 \end{equation*}
 for every $\lambda >0$, $0\leq \mu\leq c\lambda^{1-\delta}$.\\ 
  Moreover we get that for every $\delta\in[0, 1]$  and $C_0>0$,  there exists $C>0$ such that for every $\theta\in[0, \frac{d}{k}]$  there exists $C>0$ such that
 \beq\label{scestp}
 \vert\pi_{H_{nor}}(\lambda+\mu; x,x) - \pi_{H_{nor}}(\lambda; x,x)\vert \leq  C\lambda^{\frac{k(d+\theta)}{k+1} -\delta}\<x\>^{-k\theta}
 \eeq
 for $\vert\mu\vert \leq C_0\lambda^{1-\delta}$, $\lambda \geq 1$.\\
 It is enough to prove (\ref{scestp}) for $\vert x\vert \leq C\lambda^{\frac{1}{k+1}}$ with  $C>0$ large enough,
  using that  $\pi_{H}(\lambda; x,x)=\pi_{H_{nor}}(\lambda^{\frac{2k}{k+1}}; x,x)$
  and that for $\<x\>^{2k} \gtrsim\lambda$, $\pi_{H}(\lambda; x,x)$ is fast decreasing  in $x$ and $\lambda$.
  For $\vert x\vert \leq C\lambda^{\frac{1}{k+1}}$,  the bound   (\ref{scestp})  is clear.
 
 Finally,  the inequality  (\ref{splp}) is a consequence of an interpolation inequality and (\ref{scestp}).
     \end{proof}
  
  \section{Probabilistic weighted  Sobolev  estimates}\label{Sect3}
     Our main goal here is  to prove Theorem \ref{thm1}. We follow the strategy initiated in \cite{bule}
      and worked out in~\cite{PRT1} for the harmonic weighted spaces. A key step for the extension to potentials satisfying 
      (A1),  is to obtain  suitable spectral estimates.  
      In particular the probabilistic tools (concentration of measures)   can be applied  in a general setting as we have
      explained  in \cite{PRT1} and will be used  freely here, so we shall not give 
      all the details in the proofs.

 Here we shall denote  by $h = \lambda^{-1}$ ($\lambda$ is  a typical energy of the rescaled Hamiltonian $\hat H_{nor}$).    
 Using the Weyl asymptotic formula (\ref{Weylhom})  we can infer that there exist $h_0>0$, $c>0, C>0$ such that
  \beq\label{eva}
  ch^{-d}(b_h-a_h) \leq N_h \leq Ch^{-d}(b_h-a_h),\;\; \forall  h\in]0, h_0].
  \eeq
  We are using the same notations as in the introduction for the Hamiltonian $\hat H_{nor}$.
  In particular $N_h$ is the number of eigenvalues of $\hat H_{nor}$ in $I_h=]\frac{a_h}{h}, \frac{b_h}{h}[$.\\
  We are interested in  Sobolev type estimates for $u$ in the spectral subspace
  ${\cal E}_h:= \Pi_{H_{nor}}(I_h)(L^2(\R^d))$. 

  Following the strategy of \cite{bule}, repeated in \cite{PRT1}, we  need to estimate $\frac{e_{x,h}}{N_h}$ where
$$
e_{x,h} = \pi_{H_{nor}}(\frac{b_h}{h};x,x) - \pi_{H_{nor}}(\frac{a_h}{h};x,x).
$$

    \begin{lemm}\ph  For every $\delta\in[0, 1]$ there exist      $C>0$ and 
    $h_0>0$ such that   we have   
      \beq\label{Ek2}
   e_{x,h}  \leq CN_hh^{\frac{d-k\theta}{k+1}}\<x\>^{-k\theta},\;\;\; \forall h\in]0, h_0],\;\; \;\forall     \theta \in[0, \frac{d}{k}].
  \eeq
  \end{lemm}  
 \begin{proof} 
 To avoid some misunderstanding recall here  the meaning of the parameters $\lambda$ and $h$.
 If for~$H$  the typical energy is $\lambda$ then  for $H_{nor}$ the typical energy is $\lambda_{nor}$  where we have
 $\lambda_{nor}  =\lambda^{\frac{k+1}{2k}}$. Finally we have $h=\lambda_{nor}^{-1}$. So, we  have 
              $$
        e_{x,h} =    \pi_H(\left(\frac{b_h}{h}\right)^{2k/(k+1)};x,x) - 
       \pi_H(\left(\frac{a_h}{h}\right)^{2k/(k+1)};x,x),      
       $$
    and (\ref{Ek2}) is a consequence of (\ref{eva}) and  (\ref{scestp}). 
  \end{proof}

       
       As for the harmonic oscillator we get,    using (\ref{sobg}) and interpolation, 
       for every  $u\in{\cal E}_h$,
$\theta \geq 0$, $p\geq 2$
\bea\label{detsobk}
\Vert u\Vert_{L^{\infty,k\theta/2}(\R^d)} \leq C  \left(N_hh^{\frac{d-k\theta}{k+1}}\right)^{1/2}\Vert u\Vert_{L^{2}(\R^{d})} \\
 \Vert u\Vert_{L^{p,k\theta(p/2-1)}(\R^d)} \leq C \left(N_hh^{\frac{d-k\theta}{k+1}}\right)^{\frac{1}{2}-\frac{1}{p}}\Vert u\Vert_{L^{2}(\R^{d})},\nonumber
\eea
where  $L^{r,s}(\R^d) := L^r(\R^d, \<x\>^sdx)$ for $s\in\R$, $r\geq 1$.
Notice that $N_h$ is of order $(b_h-a_h)h^{-d}$.
   
   \subsection{Upper bounds}

We can get now the following result.
\begin{theo}\ph\label{thps1}
There exist $h_0\in]0, 1]$,   $c_2>0$ and $C>0$  such that if $c_1=d(1+d/2)$ 
   we have 
\begin{equation*}
{\bf  P}_{\gamma,h}\Big[u\in{\bf S}_h,\; h^{-\frac{d-k\theta}{2(k+1)}}\Vert\<x\>^{k\theta/2}u\Vert_\infty >\Lambda\Big] \leq  Ch^{-c_1}{\rm e}^{-c_2\Lambda^2},\;\;
 \forall \Lambda >0,\; \forall h\in]0, h_0], \forall \theta\in[0, \frac{d}{k}].
\end{equation*}
Moreover we can choose every $ c_2< \ell_\theta$   and $h_0$ small enough where
$$ 
\ell_\theta  = \liminf_{h\searrow 0}\inf_{x\in\R^d}\left(\frac{N_h-1}{e_{x,h}}h^{\frac{d-k\theta}{k+1}}\<x\>^{-k\theta}\right),
$$
and $\di{\inf_{\theta\in[0,d/k]}\ell_\theta >0}$. 
\end{theo}

\begin{proof}
Similarly to \cite[Theorem 4.1]{PRT1}, we first show that there exists  $c_{2} >0$ such that for every $\theta \in[0, \frac{d}{k}]$ 
every $x\in\R^d$, and every  $\Lambda>0$ we have 
\begin{equation*}
{\bf  P}_{\gamma,h}\Big[u\in{\bf S}_h,\; \<x\>^{\frac{k\theta}{2}}h^{-\frac{d-k\theta}{2(k+1)}}\vert u(x)\vert >\Lambda\Big] \leq {\rm e}^{-c_2\Lambda^2},
\end{equation*}
then a partition of unity argument yields the result. The fact that $\ell_{\theta}>0$ follows from \eqref{Ek2}.
\end{proof}

So we  get probabilistic Sobolev  estimates  improving the deterministic one
(\ref{detsobk}) with probability close to one as $h\rightarrow 0$. The improvement is "almost" of order
$N_h^{1/2} \approx\left((b_h-a_h)h^{-d}\right)^{1/2}$. Choosing $\Lambda = \sqrt{-K\log h}$ we get the next result

\begin{coro}\ph\label{kpest}
For every $K>{c_1}/{c_2}$  we have
\begin{equation*}
{\bf  P}_{\gamma,h}\left[u\in{\bf S}_h,\; \Vert u\Vert_{L^{\infty,k\theta/2}(\R^d)} >
\vert K\log h\vert^{1/2}h^{\frac{d-k\theta}{2(k+1)}}\right] \leq h^{Kc_2-c_1},\; 
\;\; \forall h\in]0, h_0], \;\;\forall \theta\in[0, \frac{d}{k}].
\end{equation*}
and
 \begin{equation*}
{\bf  P}_{\gamma,h}\left[u\in{\bf S}_h:\; \Vert u\Vert_{{\cal W}_k^{s,\infty}(\R^d)} >
 Kh^{\frac{d}{2(k+1)}-s}\vert\log h\vert^{1/2}\right] \leq h^{Kc_2-c_1},\quad 
  \forall h\in]0, h_0],\;\forall s\geq 0.
\end{equation*}
\end{coro}
Using the same argument that in \cite[Section 3.1]{PRT1}  we can deduce from  the last corollary a global probabilistic Sobolev estimate.  Let
$\theta$ a $C^\infty$ real function on $\R$ such that $\theta(t)=0$ for $t\leq a$,
$\theta(t)=1$ for $t\leq b/2$ with $0 <a < b/2$.  Define
$\psi_{-1}(t) = 1-\theta(t)$, $\psi_j(t) = \theta(h_jt) - \theta(h_{j+1}t)$ for $j\in\N$.
For every  distribution $u\in{\cal S}'(\R^d)$ we introduce  the following   dyadic Littlewood-Paley decomposition ($h_j=2^{-j}$),
\begin{equation*}
	u = \sum_{j\geq -1}u_j,\quad {\rm  with}\quad u_j = \sum_{\ell \in\N}\psi_j(\lambda_\ell)\<u,\varphi_\ell\>\varphi_\ell
\end{equation*}
and we have  $u_j\in{\cal E}_{h_j}$. \\
The following scale of spaces contains the Besov spaces $B^s_{2,\infty}$ for every $s>0$ and 
$m\geq 0$.
$$
{\cal G}^m = \big\{u\in \mathcal{S}'(\R^{d}):\; \sum_{j\geq 1}j^{m}\Vert u_j\Vert_{L^2(\R^d)} <+\infty\big\},\quad m\geq 0.
$$
We introduce probabilities $\mu_\gamma^{m}$ on  ${\cal G}^m$ (see  \cite{PRT1}  for details) as the law of
$\di{v_\gamma(\omega) = \sum_{j\geq 0}\gamma_jX_j(\omega)\varphi_j}$.
In particular we can manage such that $\mu_\gamma^{1/2}({\cal G}^{1/2}) = 1$ and
$\mu_\gamma^{1/2}({\cal G}^{m}) = 0$  for $m>1/2$.
\begin{coro}
There exists $c>0$ such that for  every $\kappa>0$ large enough, there exists a Borel subset $B^\kappa$ of
${\cal G}^{1/2}$ with $\mu_\gamma^{1/2}(B^\kappa)\geq 1-{\rm e}^{-c\kappa^2}$,  such that for every $u\in B^\kappa$ we have
\begin{equation*}
\Vert u\Vert_{\W_k^{\frac{d}{2(k+1)},\infty}(\R^d)} \leq \kappa \Vert u\Vert_{{\cal G}^{1/2}}.
\end{equation*}
\end{coro}
 Here, when $k\to \infty$, we formally recover the estimates proved by Burq-Lebeau \cite[Corollaire 1.1]{bule} in case of a compact manifold.

\subsection{Lower bounds}
Here we suppose that the random variables follow the Gaussian law ${\cal N}_\C(0,1)$ and that $\delta < 2/3$.

Now we are interested to get a lower bound for $\Vert u\Vert_{L^{\infty,k\theta/2}(\R^d)}$
 and $\Vert u\Vert_{\W_k^{s,\infty}}$. 
As we have seen in~\cite{PRT1}, a first step is to get  two sides weighted  $L^r$ estimates  for large $r$.\ligne

Our goal now is to prove the following probabilistic  improvement  of (\ref{detsobk}) when $\theta=d/k$
\begin{theo}\ph\label{2ESTr}
 There exist $0<K_0 < K_1$, $K>0$,  $c_1>0$ , $h_0>0$ such that
\begin{equation*} 
{\bf P}_{\gamma,h}\Big[u\in {\bf S}_h : \Big| \Vert u\Vert_{L^{r, \frac{d}{k}(\frac{r}2-1)}}  - {\cal M}_{r}\Big| >\Lambda\Big]
\leq 2\exp\big(-c_2N_h^{2/r}\Lambda^2\big),
\end{equation*}
for all $r\in [2, K\vert\log h\vert]$ and  $h\in]0, h_0]$.

Therefore   for every $\kappa\in]0, 1[$,  $K>0$, there exist $0<C_0 < C_1$,  $c_1>0$ , $h_0>0$ such that for all 
$r\in [2, K\vert\log h\vert^\kappa]$,  $h\in]0, h_0]$ and $\Lambda >0$  we have 
\begin{equation*}
{\bf P}_{\gamma,h}\Big[u\in{\bf S}_h : C_0\sqrt r \leq \Vert u\Vert_{L^{r, \frac{d}{k}(\frac{r}2-1)}}
 \leq C_1\sqrt r\Big] \geq 1- {\rm e}^{-c_1\vert\log h\vert^{1-\kappa}}.
 \end{equation*}
\end{theo}

Thus for $\theta=d/k$ we get a two sides weighted $L^\infty$ estimate showing that  Theorem \ref{thps1} and its corollary are sharp.
\begin{coro} 
After  a slight modification of the constants in Theorem \ref{2ESTr}, if necessary,  we get
 \begin{equation*}
 {\bf P}_{\gamma,h}\left[u\in{\bf S}_h : K_0\vert\log h\vert^{1/2} \leq \Vert u\Vert_{L^{\infty, d/2}} \leq K_1\vert\log h\vert^{1/2} \right] \geq 1- h^{c_1}
 \end{equation*}
for all  $h\in]0, h_0]$.
\end{coro}
  The key ingredient of the proof is a two side estimate of a weighted norm of the spectral function. More precisely, we have the following result, which will be proven in Appendix~\ref{A.4}.
  \begin{lemm}\ph\label{2estspf}
  There exist $0<C_0 <C_1$ and $h_0>0$ such that
\begin{equation*}  
  C_0N_hh^{\beta_{2p,\theta}}\leq \left(\int_{\R^d}\<x\>^{k\theta(p-1)}{\rm e}^{p}_{x,h}dx\right)^{1/p}
  \leq   C_1N_hh^{\beta_{2p,\theta}}, 
\end{equation*}
  for  every $p\in [1, \infty[$ and $h\in ]0, h_0]$ where $\beta_{r,\theta}=\frac{d-k\theta}{k+1}(1-\frac{2}{r}).$
  \end{lemm}
Set  $F_r(u) = \Vert u\Vert_{L^{r, \theta(r/2-1)}}$ and denote by ${\cal M}_{r}$ its median,
and by ${\cal A}_r^r = {\bf E}_h(F_r^r)$ the moment of  order~$r$. From Lemma \ref{2estspf} we get the estimates (see the proof of \cite[Theorem 4.7]{PRT1})
\begin{equation*} 
  C_0\sqrt rh^{\beta_{r,\theta}} \leq {\cal M}_{r}, {\cal A}_{r} \leq C_1\sqrt rh^{\beta_{r,\theta}} ,\;\;\;\forall h\in ]0, h_0], \; \;r\in[2, K\log N_h] . 
 \end{equation*}
 In particular, for $\theta =\frac{d}{k}$  we get 
\begin{equation*}  
   C_0\sqrt r \leq {\cal M}_{r},  {\cal A}_{r} \leq C_1\sqrt r,\;\;\;\forall h\in ]0, h_0], \; \;r\in[2, K\log N_h] ,
  \end{equation*}
  so Theorem \ref{2ESTr}  and its corollary follow from concentration of measures properties.
  
  \begin{rema}
 For the median ${\cal M}_{\infty}$ and the mean ${\cal A}_{\infty}$ of
  $F_\infty(u):=\Vert u\Vert_{L^{\infty,d/2}}$ we get the estimates
  $$
  K_0\vert\log h\vert^{1/2} \leq {\cal M}_{\infty}, {\cal A}_{\infty} \leq K_1\vert\log h\vert^{1/2},\;\;\forall h\in]0, h_0].
  $$ 
 \end{rema}
 We now  state analogous results for the norm in the Sobolev scale ${\cal W}_k^{s,r}$  for which  Theorem \ref{thm1} is a particular case. Denote by ${\cal M}_{r,s}$ a median of $\Vert u\Vert_{\cal W_k^{s,r}(\R^d)}$ and set 
 $
\dis \beta_s = \frac{d}{2(k+1)}(1-\frac{2}{r})-s.
$
Then  
  \begin{theo}\ph 
 Let $s\geq 0$. There exist $0<C_0 < C_1$, $K>0$,  $c_1>0$ , $h_0>0$ such that for all $r\in [2, K\vert\log h\vert]$ and  $h\in]0, h_0]$
 \begin{equation*} 
{\bf P}_{\gamma,h}\Big[u\in {\bf S}_h : \Big| \Vert u\Vert_{\cal W_k^{s,r}(\R^d)}  - {\cal M}_{r,s}\Big| >\Lambda\Big]
\leq 2\exp\big(-c_2N_h^{2/r}h^{-2\beta_s}\Lambda^2\big)
\end{equation*}
and 
\begin{equation*}
 C_0\sqrt rh^{\beta_s} \leq {\cal M}_{r,s} \leq C_1\sqrt rh^{\beta_s},\quad
  \forall r\in[2, K\log N].
  \end{equation*}
  
  In particular, for every $\kappa\in]0, 1[$, $K>0$ , there exist $C_0 >0$, $C_1>0$, $c_1>0$ such that 
  for  every $r\in [2, K\vert\log h\vert^\kappa]$ we have 
\begin{equation*}
{\bf P}_{\gamma,h}\left[u\in{\bf S}_h : C_0\sqrt r h^{\frac{d}{2(k+1)}(1-\frac{2}{r})}h^{-s}\leq \Vert u\Vert_{{\cal W}_k^{s,r}(\R^d)}
 \leq C_1\sqrt r h^{\frac{d}{2(k+1)}(1-\frac{2}{r})}h^{-s}\right] \geq 1- {\rm e}^{-c_1\vert\log h\vert^{1-\kappa}},
\end{equation*}
and for $r=+\infty$ we have for all  $h\in]0, h_0]$
\begin{equation*}
{\bf P}_{\gamma,h}\left[u\in{\bf S}_h : C_0\vert\log h\vert^{1/2}h^{\frac{d}{2(k+1)}-s} \leq\Vert u\Vert_{{\cal W}_k^{s,\infty}(\R^d)}
\leq C_1\vert\log h\vert^{1/2}h^{\frac{d}{2(k+1)}-s} \right] \geq 1- h^{c_1}.
\end{equation*}
 \end{theo}

 \begin{rema}
 For more general probability laws satisfying  the  Gaussian concentration estimate
 we can  get a weaker lower bound  as in \cite[Theorem 4.13]{PRT1}.
 
 \end{rema}

    \subsection{No random smoothing in $L^2$-scale spaces}
    This  subsection is  a generalization to non harmonic potentials of results proved by Poiret \cite[Section 4]{poiret}, and we refer to this paper for more details.
    
    Let $\di{v_\gamma = \sum_{j\geq 0}\gamma_j\varphi_j}$  be such that $v_\gamma\in{\cal S}^\prime(\R^d)$.
    Denote by $\dis v_\gamma^\omega = \sum_{j\geq 0}\gamma_jX_j(\omega)\varphi_j$ a randomization
    of $v_\gamma$ with $\nu$ satisfying (\ref{Blip}).\\
    Recall that ${\cal H}_k^s = Dom(-\Delta + \vert x\vert^{2k})^{\frac{k+1}{2k}s}$ and ${H}^{\s} = Dom(-\Delta  )^{\frac{\s}{2}}$, then 
    \begin{theo}\ph\label{nors}
    If for $s\geq0$, $v_\gamma\notin{\cal H}_{k}^{s}(\R^d)$ then
    a. s.  $v_\gamma^\omega\notin H^{\frac{k+1}{k}s}(\R^d)$   and 
      {${v_\gamma^\omega\notin L^2(\R^d,\vert x\vert^{(k+1)s}dx)}$}
      \end{theo}
    The first  step in  the proof is the following
    \begin{prop}
    For every $s\geq 0$ there exist $C_1(s)>0, C_2(s)>0$ such that
    \beq\label{eqrep1}
    C_1(s)\lambda_n^s \leq \Vert\Delta^s\varphi_n\Vert_{L^2(\R^d)} \leq C_2(s)\lambda_n^s,\; \forall n\geq 0,
    \eeq
    \beq\label{eqrep3}
    C_1(s)\lambda_n^{s/(2k)} \leq \Vert\vert x\vert^s\varphi_n\Vert_{L^2(\R^d)} \leq C_2(s)\lambda_n^{s/(2k)},\; \forall n\geq 0.
    \eeq
    \end{prop}
    \begin{proof}
    By  a standard scaling argument, estimate (\ref{eqrep1}) is equivalent to the semiclassical estimate
    \beq\label{eqrep2}
    C_1(s)\leq \Vert(h\Delta)^s\psi_h\Vert_{L^2(\R^d)} \leq C_2(s),\; \forall h\in]0, 1],
    \eeq
    where $(-\h^2\Delta +V)\psi_h =\lambda_h\psi_h$ and $\di{\lim_{h\rightarrow 0}\lambda_h =1}$.\\
    In (\ref{eqrep1}) or (\ref{eqrep2})  the upper bound is obvious and we only need to prove the lower bound. \\
    Assume that the lower bound in \eqref{eqrep2}  is not satisfied. 
    Let $\mu_{sc}$ be the semiclassical measure of the semiclassical state $\psi_h$. It is well known that the support of $\mu_{sc}$
    is in the energy surface $H^{-1}(1)$ and is invariant by the classical Hamiltonian flow $\Phi_H^t$ of $H$. Denote by $A (x,\xi)=\vert\xi\vert^{2s}$. Then for some sequence of $h$  we have
    \begin{equation}\label{lim}
      \lim_{h\rightarrow 0}\Vert\hat A(h)\psi_h\Vert_{L^2(\R^d)} =0.
    \end{equation}
    By definition of the semiclassical measure $\mu_{sc}$, for some subsequence $h_{n}$ we have
    \begin{equation*}
    \<\hat A(h_{n})\psi_{h_{n}}, \psi_{h_{n}}\>\longrightarrow \int_{\R^{2d}}A(x,\xi)\text{d}\mu_{sc}(x,\xi),
    \end{equation*}
    therefore, from \eqref{lim}, using that $A\geq 0$ we get that the support of $\mu_{sc}$ is in $A^{-1}(0)\cap H^{-1}(1)$.
    So if $X_0$ is in supp($\mu_{sc}$) then for every $t\in\R$, $\Phi_H^tX_0\in {\rm supp(}\mu_{sc}$).
     We get a contradiction  because $A^{-1}(0)$ cannot contain a global  classical trajectory.
     The same argument applies to the observable $A(x,\xi) = \vert x\vert^s$, giving \eqref{eqrep3}.
    \end{proof}
    
        Now we prove Theorem \ref{nors} following \cite{poiret}.
        
    \begin{proof}[Proof of Theorem \ref{nors}]
  Let $\chi$ a smooth cutoff, $0\leq \chi\leq 1$, ${\rm supp}(\chi)\subseteq [-1, 2]$, $\chi=1$ on $[0, 1]$.
    Introduce  the notations:
    \begin{eqnarray*}
      \sigma_N^2 &=&    \E\Big[   \big\Vert\chi\Big(\frac{H}{N}\Big)v_\gamma^\omega\big\Vert_{\H^{s}(\R^d)}^2\Big]= \sum_{n\geq 1}\chi^2\Big(\frac{\lambda_n}{N}\Big)\vert\gamma_n\vert^2\lambda_n^{s}\\
    s_N(\omega) &= &\big\Vert\chi\Big(\frac{H}{N}\Big)v_\gamma^\omega\big\Vert_{H^{s}(\R^d)} \\
    M(\omega) &= & \sup_{N\geq 1}s_N(\omega).
    \end{eqnarray*}
  The event   $[M=+\infty]$ is a tail event for the sequence of independent random variables $\{X_n\}_{n\geq 0}$. Thus, according to the (0-1) law,  to prove  that $\P[M=+\infty]=1$ it is enough to prove that $\P[M=+\infty] >0$.\\
    We will use the following Paley-Zygmund inequality: if $X\geq 0$  is a random variable and $\lambda\geq 0$, 
     \begin{equation*}
    \P\Big[X\geq \lambda E(X)\Big] \geq (1-\lambda)^2\frac{(EX)^2}{EX^2}.
   \end{equation*}
    Using \eqref{eqrep1}  and properties of the sequence $\{X_n\}_{n\geq 0}$  we get easily, for  some $c>0, C>0$, 
$$
    \E\big[\,s_N^2\,\big] \geq c\sigma_N^2\;\;\, {\rm and}\;\;  \,    \E\big[\,s_N^4\,\big] \leq C\sigma_N^4,\;\; \forall N\geq 1.
      $$
    Hence using the Paley-Zygmund inequality, we get
    $$
    \P\left[M^2 \geq \frac{c\sigma_N^2}{2}\right]  \geq \frac{c^2}{4C^2},\;\;\forall N\geq 1.
    $$
    The result follows from the monotone convergence theorem.
    \end{proof}
 
    \section{Random quantum ergodicity}\label{Sect4}
   We adapt here in the  context of semi-classical Schr\"odinger operators on $\R^d$   results  proved  
     for the Laplace operator on compact manifolds \cite[Theorem 3]{bule}.
     
   The spectral results which are necessary in the proof come from  \cite{hmr}, and have already been used there to obtain   quantum ergodicity results. For an introduction to quantum ergodicity, we refer to Zworski \cite[Chapter 15]{Zworski}.
    \subsection{Spectral preparations}
    Let   us introduce  now the assumptions  for this section. Let   $d\geq 2$ and let $J_h=hI_h= ]a_h, b_h]$
    be such that 
    \begin{equation}\label{length}
      \lim_{h\rightarrow 0}a_h=  \lim_{h\rightarrow 0}b_h=\eta\qquad \text{and}\qquad  
    \lim_{h\rightarrow 0}\frac{b_h -a_h}{h}=+\infty.
    \end{equation}
    For a smooth symbol $A(x,\xi)$ in $\R^{2d}$, we denote by $\hat A(h)$ its  Weyl $h$-quantization    defined by 
   \begin{equation*}
   \hat{A}(h)f(x)=\hat{A}(x,hD)f(x)=\frac{1}{(2\pi h)^{d}}\int_{\R^{d}}\int_{\R^{d}}e^{i(x-y)\xi/h}A(\frac{x+y}{2},\xi)f(y)dyd\xi,
   \end{equation*}
then if $A$ is real, the operator   $\hat A(h)$ is self-adjoint in $L^2(\R^d)$.
   
     Assume that  $H(x,\xi)=|\xi|^{2}+V(x)$, with $(x,\xi)\in \R^{2d}$ where $V$ satisfies the assumptions (A1). Then from~\cite{hmr} (and its bibliography) it is known that:
    \begin{enumerate}[(i)]
\item In $J_h$ the spectrum of $\hat H(h)$ is discrete.
\item   The number $M_h$ of eigenvalues of $\hat H(h)$ in $J_h$ satisfies the Weyl-H\"ormander law:
   $$
   M_h =(2\pi h)^{-d}\int_{[(x,\xi); H(x,\xi)\in J_h]}dxd\xi + {\cal O}(h^{1-d}).
$$
Therefore, thanks to Assumption $(A1)$ on $V$ we can make the change of variables $y=h^{1/(k+1)}x$, $\eta=h^{k/(k+1)}\xi$ and get, thanks to \eqref{eqN}, that $N_{h}\sim M_{h}$ when $h\longrightarrow 0$.
\item Denote by $S(1)$ the set of the symbols
  $$
   S(1) =
   \big\{A\in C^\infty(\R^{2d}),\;  \forall\alpha, \forall\beta,\;\; \sup_{(x,\xi)}\vert\partial_x^\alpha\partial_\xi^\beta A(x,\xi)\vert <+\infty\big\}.
   $$
By the Calderon-Vaillancourt theorem, there exist $c>0$ and $\nu_{d}\in \N$ such that for all $A\in S(1)$
\begin{equation*}
\|\hat{A}\|_{L^{2}\to L^{2}}\leq c \|A\|_{W^{\nu_{d},\infty}(\R^{2d})}.
\end{equation*}
Therefore, for later purpose, we introduce
\begin{equation}\label{defsb}
   S_{b}(1) = \big\{A\in S(1),\;    \|A\|_{W^{\nu_{d},\infty}(\R^{2d})}\leq 1\big\}.
\end{equation}
  Next, let $\Pi_h$  be the spectral projector of $\hat H(h)$ on $J_h$ and $A\in S(1)$, then we have
  \begin{equation*}
   {\rm Tr}\big(\hat A(h)\Pi_h\big) =  (2\pi h)^{-d}\int_{[(x,\xi); H(x,\xi)\in J_h]}A(x,\xi) dxd\xi + {\cal O}(h^{1-d}).
\end{equation*}
\end{enumerate}
   Denote by  $dL_\eta$  the Liouville  probability measure on  the energy  surface $H^{-1}(\eta)$, and recall  that 
   \begin{equation*}
      dL_\eta = C_\eta\frac{d\Sigma_\eta}{\vert\nabla H\vert},
   \end{equation*}
   where  $ d\Sigma_\eta$ is the Euclidean measure on $\Sigma_\eta:=H^{-1}(\eta)$
     and $C_\eta$ the normalization constant. Then we can prove that 
     \begin{equation*} 
      \frac{ {\rm Tr}\big(\hat A(h)\Pi_h\big) }{M_h} = L_\eta(A) + {\cal O}\Big(\frac{h^{1-d}}{M_h}+b_h-a_h\Big)= L_\eta(A)  
    + {\cal O}\Big(\frac{h}{b_h-a_h} +b_h-a_h\Big),
     \end{equation*}
     so that with our assumptions we have
     \begin{equation}\label{aver2}
        \lim_{h\rightarrow 0}\frac{ {\rm Tr}\big(\hat A(h)\Pi_h\big) }{M_h} = L_\eta(A).
     \end{equation}

     As above,  ${\cal E}_h$ is the range of $\Pi_h$ and  ${\bf S}_h$ is the (complex) unit sphere of ${\cal E}_h$, which is of dimension $M_{h}$. In the particular case where  ${\bf Q}^{(u)}_h$ is  the uniform probability on ${\bf S}_h$, using that the canonical probability on the sphere is invariant under the unitary group, we can see easily that the corresponding expectation satisfies 
\begin{equation*}      
      {\bf F}^{(u)}_h\Big[\<v, \hat A(h) v\>\Big] = \frac{ {\rm Tr}\big(\hat A(h)\Pi_h\big) }{M_h},\;\; {\rm for \; every}\;\;   v\in{\bf S}_h.
     \end{equation*}
      We shall see later (in \eqref{trace}) that this is  still true, up to  a small  error,  for more general probabilities.

      \subsection{Proof of Theorem \ref{rqer} } Recall the definition \eqref{defpsi}
      \begin{equation}\label{defpsi2}
         \psi^h_j(x) = h^{-\frac{d}{2(k+1)}}\varphi_j(h^{-\frac{1}{k+1}}x):=\mathscr{L}_h\, \phi_{j}(x).
      \end{equation}
We   define the random vector  in ${\cal E}_h$ 
\begin{equation*} 
\tau(\omega) : = \frac1{\sqrt{M_{h}}}\sum_{j\in\Lambda_{h}} X_j(\omega)\psi^{(h)}_{j},
 \end{equation*}
 and we assume  that the law of the r.v. $X=(X_{j})_{j\geq 1}$ satisfy  \eqref{Blip}. We then  consider  ${\bf Q}_h$  the probability on ${\bf S}_h$   defined   by
 \begin{equation*} 
  \int_{{\bf S}_h}f(u)\text{d}{\bf Q}_{h}(u) =  \int_\Omega f\left(\frac{\tau(\omega)}{\Vert \tau(\omega)\Vert_{L^{2}(\R^{d})}}\right)\text{d}\P(\omega),
 \end{equation*}
for all measurable and bounded function $f:{\bf S}_h\longrightarrow \R$. In the sequel we denote by ${\bf F}_{h}$ the expectation associated to this probability.\ligne
     
     Theorem \ref{rqer} will be implied by the  following Proposition, with the unitary transformation $\mathscr{L}_{h}$. Recall the definition \eqref{defsb} of $S_{b}(1)$, then

    \begin{prop}\ph\label{scquerg} 
Under the previous assumptions there exist $c,C>0$ so that for all $r\geq 1$, $h\in]0, 1]$  and $A\in S_{b}(1)$
  \begin{equation}\label{rqest}
  {\bf Q}_{h}\Big[ \,u\in {\bf S}_{h} : |\<u, \hat A(h) u\>-L_{\eta}(A)|>r   \Big]\leq C\e^{-cM_hr^{2}}.
  \end{equation}    
    As a consequence, there exists $C>0$ so that for all $p\geq 2$, $h\in ]0, 1]$
     and $A\in S(1)$
    \begin{equation*}
    \big\|\<u, \hat A(h) u\>-L_{\eta}(A)\big\|_{L^{p}_{{\bf Q}_{h}}}\leq CM_h^{-1/2}\sqrt{p}.
    \end{equation*} 
      \end{prop}
      
      \begin{rema}
      Actually, this result holds true under  more general hypotheses on $H$. Namely, it is enough to assume that $H$ is a real smooth symbol on $\R^{2d}$ so that 
          \begin{equation*}
      \left\{
      \begin{aligned}
         &\bullet \;\;  \vert\partial_x^\alpha\partial_\xi^\beta H(x,\xi)\vert \leq C_{\alpha, \beta},\;\; \vert\alpha+\beta\vert \geq m  \;\;\text{if}\;\;  m  \;\;\text{is large enough}  \\
 &\bullet \;\;H^{-1}[\eta-\varepsilon_0, \eta+\varepsilon]\;\; {\rm is\; compact\; in}\; \;\R^{2d}   \\
&\bullet \;\;\eta \; {\rm is\; non \;critical\; for}\; H.
\end{aligned}
\right.
    \end{equation*}
        \end{rema}
      
               \begin{proof}
  Let $A\in S_{b}(1)$. To begin with, we can assume that $\hat A(h)\geq 0$, since we can consider the operator $\hat A(h)+C$, where $C>0$ is some large constant. Then  observe that for all $u,v\in {\bf S}_{h}$, 
      \begin{equation*}
      \big|\<u, \hat A(h) u\>-\<v, \hat A(h) v\>\big|\leq 2\|\hat{A}(h)\|_{L^{2}\to L^2}\|u-v\|_{L^{2}(\R^{d})}\leq c \|u-v\|_{L^{2}(\R^{d})},
      \end{equation*}
      thus we can apply \cite[Proposition 2.12]{PRT1}, which gives $K,\kappa>0$ so that for all $r>0$ and $h\in]0,1]$
      \beq\label{supconcent}
{\bf Q}_{h}\Big[ \,u\in {\bf S}_h : \vert  \<u, \hat A(h) u\>-\mathcal{M}\vert > r\,\Big] 
\leq K{\rm e}^{-\kappa M_hr^2},
\eeq
where $\mathcal{M}$ is a median for $u\longmapsto \<u, \hat A(h) u\>$. Next, by \eqref{supconcent}
\begin{eqnarray*}
\big|\,{\bf F}_{h}\big[\, \<u, \hat A(h) u\>\,\big]-\mathcal{M}\,\big|&\leq &{\bf F}_{h} \big[\,| \<u, \hat A(h)u\>-\mathcal{M}  |  \,\big]\nonumber \\
&=&\int_{0}^{+\infty} {\bf Q}_{h}\Big[ \,u\in {\bf S}_h : \vert  \<u, \hat A(h) u\>-\mathcal{M}\vert > r\,\Big] \text{d}r\leq CM_h^{-1/2}.
\end{eqnarray*}
In order to complete the proof, by \eqref{aver2} and the previous lines, it is enough to show that 
 \begin{equation}\label{trace}
       {\bf F}_h\big[\,\<u, \hat A(h) u\>\,\big] - \frac{{\rm Tr}(\hat A(h)\Pi_h)}{M_h} \longrightarrow 0,
 \end{equation}
when $h\longrightarrow 0$. Here the random vector $\tau$ which defines ${\bf Q}_{h}$ reads 
\begin{equation*}
\tau(\om)=\frac1{\sqrt{M_h}}\sum_{j\in \Lambda_h}X_{j}(\om)\psi^{(h)}_{j},
\end{equation*}
so that 
\begin{equation*}
\|\tau(\om)\|^{2}:=\|\tau(\om)\|^{2}_{L^{2}(\R^d)}=\frac{1}{M_h}\sum_{j=1}^{M_h}|X_{j}(\om)|^{2}.
\end{equation*}
By definition 
\begin{equation*}
 {\bf F}_h\big[\,\<u, \hat A(h) u\>\,\big] =\int_{{\bf S}_{h}}\<v,\hat A(h) v\>\text{d}{\bf Q}_{h}(v)=\int_{\Omega}\<\frac{\tau(\om)}{\|\tau(\om)\|},\hat A(h) \frac{\tau(\om)}{\|\tau(\om)\|}\>\text{d} \P(\omega).
\end{equation*}
Denote by 
\begin{equation*}
\Omega_{h}=\big\{\,\om \in \Omega : \big| \|\tau(\om)\|^{2}_{L^{2}}-1  \big|\leq M_h^{-1/2}\,\big\}.
\end{equation*}
By \cite[Lemma 2.11]{PRT1}, $\P(\Omega^{c}_{h})\leq C\e^{-cM_h^{1/2}}$. We write
\begin{multline}\label{esp}
\int_{\Omega}\<\frac{\tau}{\|\tau\|},\hat A(h) \frac{\tau}{\|\tau\|}\>\text{d} \P(\omega)-\int_{\Omega}\<\tau,\hat A(h) \tau\>\text{d} \P(\omega)=\\
\int_{\Omega_{h}}\frac{1-\|\tau\|^{2}}{\|\tau\|^{2}}\<\tau,\hat A(h) \tau\>\text{d} \P(\omega)+\int_{\Omega^{c}_{h}}\frac{1-\|\tau\|^{2}}{\|\tau\|^{2}}\<\tau,\hat A(h) \tau\>\text{d} \P(\omega):=C_{1,h}+C_{2,h}.
\end{multline}
Firstly,  by definition of $\Omega_{h}$ we get $|C_{1,h}|\leq M_h^{-1/2}\|\hat A(h)\|_{L^{2}}$. Secondly, by Cauchy-Schwarz
\begin{equation*}
|C_{2,h}|\leq \|\hat A(h)\|_{L^{2}}\int_{\Omega^{c}_{N}} (1+\|\tau\|^{2})\text{d} \P(\omega)\leq \|\hat A(h)\|_{L^{2}}\P^{1/2}(\Omega^{c}_{h})\big(\int_{\Omega}(1+\|\tau\|^{2})^{2}\text{d} \P(\omega)\big)^{1/2}.
\end{equation*}
It remains to check that $\dis \int_{\Omega}\|\tau\|^{4}\text{d} \P(\omega) \leq C$. By \eqref{quad}, there exists $\eps>0$ so that $\dis  \int_{\Omega}\e^{{2\eps}X^{2}_{j}(\om)}\text{d} \P(\omega)=C_{\eps}<+\infty$. Thus, by the inequality   $\eps^{2} |\tau|^{4}\leq \e^{2\eps |\tau|^{2}}$ and by Jensen 
\begin{equation*}
\eps^{2}\int_{\Omega}\|\tau\|^{4}\text{d} \P(\omega) \leq \int_{\Omega}\e^{\frac{2\eps}{M_h}\sum_{j=1}^{M_h}X^{2}_{j}(\om)}\text{d} \P(\omega) \leq \big(\int_{\Omega}\e^{2\eps\sum_{j=1}^{M_h}X^{2}_{j}(\om)}\text{d} \P(\omega) \big)^{1/M_h}=C_{\eps}.
\end{equation*}
Therefore $C_{1,h}+C_{2,h}\longrightarrow 0$.

Finally, observe that 
\begin{equation*}
\int_{\Omega}\<\tau,\hat A(h) \tau\>\text{d} \P(\omega)=\frac1{M_h}\sum_{j=1}^{M_h}\<\psi^{(h)}_{j},\hat A(h) \psi^{(h)}_{j}\>= \frac{{\rm Tr}(\hat A(h)\Pi_h)}{M_h},
\end{equation*}
and thanks to \eqref{esp}, we get \eqref{trace}.
      \end{proof}
       
      Now we can prove Theorem \ref{rqer}.  To do that we apply Proposition  \ref{scquerg}
      to the semiclassical Hamiltonian $\hat H_\eps(h)=-h^2\Delta +V_\eps$ depending
      smoothly on $\eps\in[0, 1]$. 
      This reduction is  allowed using the following Lemma
      \begin{lemm}\ph\label{clue} Consider the $(\psi^{h}_{j})$ defined by \eqref{defpsi2}. Then for every $A\in S(1,k)\subset S(1)$ we have 
           $$
           \<\varphi_j, \hat A\varphi_{\ell}\> = \<\psi^h_j, \hat A(h)\psi^h_{\ell}\> + {\cal O}(h^\infty)
           $$
        uniformly in $j,{\ell}\in\Lambda_h$.
      \end{lemm}
      \begin{proof}
      Let us introduce the unitary transformation
      $\mathscr{L}_hu(x) = h^{-\frac{d}{2(k+1)}}u\left(h^{-\frac{1}{k+1}}x\right)$. We have
      $$
      \<\varphi_j, \hat A\varphi_{\ell}\> = \<\psi^h_j, \mathscr{L}_h\hat A\mathscr{L}_h^{-1}\psi^h_{\ell}\>.
      $$
      The  1-Weyl symbol  $A_h$ of  $\mathscr{L}_h\hat A\mathscr{L}_h^{-1}$ is
      $A_h(x,\xi) = A\left(h^{-\frac{1}{k+1}}x, h^{\frac{1}{k+1}}\xi\right)$.  So its $h$-Weyl symbol 
      is $A(h)(x,\xi) =  A\left(h^{-\frac{1}{k+1}}x, h^{-\frac{k}{k+1}}\xi\right)$.
      Now using that $A$ is quasi-homogeneous outside $(0,0)$, we have for every $h\in]0, 1]$ and
      $\vert x\vert +\vert\xi\vert \geq \eps_0$,  that $A(h)(x,\xi) =A(x,\xi)$.
      
      It is known that, for every $j\in\Lambda_h$,  the semi-classical wave front set  (or the frequency set, see for example \cite{ro2}) of $\psi^h_j$ 
      is in a small neighborhood of $H_0^{-1}(\tau)$ which  do not contain a neighborhood of $(0,0)$ for 
      $\tau$ large enough. Then the  Lemma is proved.
         \end{proof}
      
      Now  we easily get  the following result, using the Borel-Cantelli Lemma.\\
           Assume that the hypothesis of Theorem \ref{rqer} are satisfied.\\  
      Let $\{h_j\}_{j\geq 0}$ be  a sequence of positive real numbers converging to 0 as $j\rightarrow +\infty$.
       Define  the compact metric space $X =\prod_{j\in\N}{\bf S}_{h_j}$  equipped with the product probability
       $$\di{{\bf P} = \otimes_{j\in\N}{\bf P}_{h_j}}.$$
      Let $u\in X$, $u=\{u_j\}_{j\in \N}$ where $u_j\in {\bf S}_{h_j}$. For any $A\in S(1,k)$,
      $u\mapsto \<u_j, \hat Au_j\>$ defines a sequence of random variables on $X$. 
       \begin{coro}\ph\label{rqe}
       Assume that $d\geq 2$ and that 
            $\di{\sum_{j\geq 0}{\rm e}^{-\eps h_j^{1-d}}<+\infty}$   for every $\eps >0$.
             Then 
   \begin{equation*}
               {\bf P}\Big[u\in X,\;\;\lim_{j\rightarrow +\infty}\<u_j, \hat Au_j\> = L_\eta(A),\;\forall A\in S(1,k)\Big] = 1.
   \end{equation*}
     \end{coro}
    \begin{proof}
    Denote by $f_h(u) = \<\Pi_hu,\hat A \Pi_hu\>$. 
Notice that the random variable $f_{h_j}$  depends only on $u_j=\Pi_{h_j}u$ so we get 
 $$
{\bf P}\Big[u=\{u_j\} \;:\; \vert f_{h_j}(u)-L_\eta(A)\vert \geq \varepsilon\Big] 
 = {\bf P}_{h_j}\Big[\vert f_{h_j} - L_\eta(A)\vert\geq \eps\Big].
$$
So applying \eqref{rqest} and the Borel-Cantelli Lemma to the independent random
variables $\{f_{h_j}\}_{j\in\N}$ we get  the conclusion.
    \end{proof}
    
      \subsection{Proof of Theorem \ref{QUE} }
      
        In (\ref{length})  the  condition  on  $b_h-a_h$ (the  length  of the  spectral  windows)
      is too restrictive to prove  the Theorem  (QUE). To achieve the  proof  we need  to  relax the condition (\ref{length}), since the harmonic  oscillator does not satisfy \eqref{length}. Actually, to isolate  its  eigenvalues  we need  to consider  spectral windows $[a_h, a_h +2h[$.\\
      To enlighten the  discussion   we  present  a more  general  setting, inspired  from  the  paper
      \cite{zel3} using   semi-classical  spectral results  proved in \cite{hr, pr}.
      
      Besides  the general assumptions  on the  classical  hamiltonian  $H$  we shall  consider  the  two
      following  particular  cases:\medskip
      \begin{enumerate} 
             \item[(Per)]\qquad  There  exists  $\epsilon_0>0$  such that in $H_0^{-1}[\eta-\epsilon_0, \eta+\epsilon_0]$
       the  hamiltonian flow $\Phi_{H_0}^t$ defined by $H_0$  is periodic with period $T=2\pi$  and 
       $T$  is a primitive period (there  exists  no  periodic  path  with period in $]0, T[$).
 \item[(APer)]\qquad    On the energy shell $\Sigma_\eta=H_0^{-1}(\eta)$  the   set of periodic  points is of
        measure 0 (for the Liouville measure).
        \end{enumerate}
    \medskip
        
       For instance,  Assumption (Per) is fulfilled   with   $H_0(x,\xi) = \frac{1}{2}(\vert \xi\vert^2 + \vert x\vert^2) $  and  (Aper)  is fulfilled
          for $d\geq 2$  
             with   $\di{H_0(x,\xi) = \frac{1}{2}(\vert\xi\vert^2 + \sum_{1\leq j\leq d} \omega_j^2x_j^2)}$,
             $\omega_j >0$, 
           if  $(\omega_1, \cdots, \omega_d)$ are  independent  on $\Z$.
      
      Denote  $\sigma(\hat H)$  the  spectrum  of $\hat H$. Remark  that  for  the  isotropic  harmonic
       oscillator  the  eigenvalues  are  very  degenerate ($d\geq 2$)  and  if  $(\omega_1, \cdots, \omega_d)$ are  independent  on $\Z$  the eigenvalues   are  non degenerate.
       \\
      If  condition (Per)  is satisfied, 
      it follows  from \cite{hr}  that there  exist  $\alpha\in\Z$, $\gamma\in\R$, $C_0\geq 0$  such that
      $$
      \sigma(\hat H)\cap]\eta-\epsilon_0, \eta+\epsilon_0[ 
      \subseteq \bigcup_{\ell\in \Z}]e_{\ell, h}-C_0h^2, e_{\ell, h}+C_0h^2]
      $$
      where $e_{\ell, h} = (\ell +\frac{\alpha}{4})h +\gamma$.\\
      For  the harmonic oscillator  we have $\alpha = 2$ (Maslov-Morse index), $\gamma=0$ and $C_{0}=0$.\\
      If  condition (Aper)  is satisfied, it follows  from \cite{pr}  that for  every $\delta>0$
      there is   an infinite  number  of  eigenvalues   of $\hat H$  in
      $J_{\delta, h} = [a_h, a_h+\delta h[$.  More  precisely  there  exists $C>0$  such that 
      $$
      \#\ \big\{\sigma(\hat H)\cap J_{\delta, h}\big\} = C\delta h^{d-1} + o(h^{d-1}).
      $$
      Our  probabilistic   results  concerning   quantum ergodicity are obtained  by combining  probability technics with     the following  spectral  semi-classical   results. 
         \begin{prop}\ph\label{prop45}
      We have  the  following  asymptotic  limit (\ref{aver2})
      $$
        \lim_{h\rightarrow 0}\frac{ {\rm Tr}\big(\hat A(h)\Pi_h\big) }{M_h} = L_\eta(A)
        $$
        for  the following  choice  of  the  spectral  window:\\
        (i)  in  the  "general" case, $J_h = [a_h, a_h + h\ell(h)[$, $\di{\lim_{h\rightarrow 0}\ell(h) = +\infty}$;\\
        (ii) under  assumption {\rm (Per)},   $J_h = ]e_{\ell, h}-C_0h^2, e_{\ell, h}+C_0h^2]$  such that
        $\di{\lim_{h\rightarrow 0, \ell\rightarrow +\infty}e_{\ell, h} = \eta}$;\\
        (iii)  under  assumption  {\rm (Aper)}, $J_h = [a_h, a_h + \delta h]$,   for  any $\delta>0$
         where $\lim_{h\rightarrow 0}a_h =\eta$.
      \end{prop}

      The previous proposition can be easily deduced from the works~\cite{hr, pr, zel3}.
      
    \begin{proof}[Proof of Theorem \ref{QUE}] We apply  Proposition~\ref{prop45} to the harmonic oscillator which satisfies the assumption $(ii)$. Then we also have the result of Proposition~\ref{scquerg}.

  Every $B\in \mathcal{B}$ can be identified with $\big\{B_j\big\}_{j\geq 1}$ where 
  $B_j\in U(N_{h_j})$. The random variables~$D_j$  are independent  and $D_j$ depends only on $B_j$. So for every $r>0$ we have
\begin{eqnarray*} 
  \rho\big[D_j(B) >r\big] = \rho_j\big[ D_j(B) >r\big] &= &
  \rho_j\Big[\,\exists \ell\in \llbracket 1,N_{h_j} \rrbracket, \;\big\vert\<\varphi_{j,\ell}, \hat A\varphi_{j,\ell}\> -L_{2}(A)\big\vert>r\Big]
  \nonumber\\
  &\leq & \sum_{1\leq \ell\leq N_{h_j}}\rho_j\Big[\,\big\vert\<\varphi_{j,\ell}, \hat A\varphi_{j,\ell}\> -L_{2}(A)\big\vert>r\Big]
 \nonumber\\
&\leq & N_{h_j}{\bf P}_{h_j}\Big[\,\big\vert\<u ,\hat A(h_j)u\> -L_{2}(A)\big\vert>r-Ch_j^M\Big].
  \end{eqnarray*}
  In the last line we have used Lemma \ref{clue} where $C\geq 0$ and $M$ arbitrary large,
   and that the probability~${\bf P}_{h}$ is the push-forward of the Haar measure  of $U(M_h)$  by the maps:  $U(N_{h})\ni M\mapsto Mv\in{\bf S}_h$,  for any  $v\in {\bf S}_h$.\\
  So using Proposition \ref{scquerg}  we get, with  positive constants  $C_1, C_2, C_3$, 
  $$
  \rho\big[D_j(B) >r\big]  \leq C_1j^{d-1}\exp\Big[-C_2j^{d-1}(r -C_3j^{-M})^2\Big].
  $$
  In  particular for any $d\geq 2$ we get 
  $$
  \sum_{j\geq 1}\rho\big[D_j(B) >r\big]  <+\infty
  $$
  and the result is a consequence of the Borel-Cantelli Lemma. 
   \end{proof}
   \begin{rema}
   Our  proof of  unique quantum  ergodicity   for random  bases can  be adapted to  prove analogous  results 
   for  the Laplace-Beltrami operator  on compact manifolds. Notice  that  our  method  is  different  from
    the method  used in \cite{zel3, map}. 
    In particular  we  do not  use  the Szeg\"o   limit  theorem like  in \cite[prop. 1.2.4]{zel3}  so that  we get  a slightly
     better  result  in the  aperiodic  case. Let us  formulate  our  result    for  non isotropic  harmonic oscillator.\\
     Let $\di{\hat H = \frac{1}{2}\big(-\triangle + \sum_{1\leq j\leq d} \omega_j^2x_j^2\big)}$, with $\omega_j>0$ and where 
      $\{\omega_1,\cdots, \omega_d\}$  are $\Z$-independent. \\
      The  eigenvalues  of $\hat H$ are
      $\di{\lambda_\alpha = \sum_{1\leq j\leq d}(\alpha_j + \frac{1}{2})\omega_j}$.
      Denote by $I_{\delta, k} = [\lambda_0 +k\delta, \lambda_0 +(k+1)\delta[$, $\delta>0$, $k\in\N$.
     For  every $\delta>0$,  the  number   of eigenvalues of $\hat H$ in  $I_{\delta, k}$ is $N_{\delta, k}\approx k^{d-1}$, $k\rightarrow +\infty$. Let  ${\cal E}_{\delta,k}$  be the linear  space spanned  by the eigenfunctions of $\hat H$ with  eigenvalues in $I_{\delta, k}$.  \\
     Denote  by ${\cal B}_\delta$  the  set  of orthonormal  bases of  $L^2(\R^d)$  obtained by choosing  an  orthonormal  basis  ${\cal B}_\delta$ in each ${\cal E}_{\delta,k}$.  Like in the proof of Theorem \ref{QUE},    ${\cal B}_\delta$   is  equipped  with  a probability  measure $\rho_\delta$. Then by
     applying Proposition~\ref{prop45}  and Proposition~\ref{scquerg}, $\rho_\delta$-almost surely,   an
     orthonormal    basis of $L^2(\R^d)$  in 
      ${\cal B}_\delta$  is QUE. \\
      Notice that every  $\psi\in  {\cal E}_{\delta,k}$  is  a $\delta$-quasimode  for $\hat H$:
      $$
      \hat H\psi = (\lambda_0+k\delta)\psi + O(\delta)\Vert\psi\Vert_{L^2(\R^d)}.
      $$

        \end{rema}
 \begin{rema}
   For Schr\"odinger operators with super-quadratic  potentials we could also get a similar result, considering orthonormal basis of quasi-modes (approximated eigenfunctions),  using Proposition\;\ref{scquerg}
    and  taking  the appropriate power  of $-\triangle + V$.
 \end{rema}
  \begin{rema}
  We have supposed that the observables $A$ are of order 0. More generally, we say that $A$ is of order $m$, $m>0$,
if   $A$ is an asymptotic sum of quasi-homogeneous symbols of degree $\leq  m$  (for more details see \cite{hero, ro1}).  We say that $A$ is  quasi-homogeneous of degree
 $m\in\R$ if 
 $A(rx, r^{k}\xi) = r^{(k+1)m}A(x,\xi)$ for $r\geq 1$, $\vert(x,\xi)\vert >\epsilon$.
 So if $A$  quasi-homogeneous   of degree $m$  then $ H_{nor}^{-m/2} A H_{nor}^{-m/2}$ is quasi-homogeneous   of degree $0$.
 So the previous results holds true for $\frac{\<\varphi_j, \hat A\varphi_j\>}{\omega_j^m}$,
 where $\omega_j = \lambda_j^{(k+1)/2k}$.

  \end{rema}


   \appendix
     \section{Some point-wise spectral estimates for confining potentials }
     The aim of this appendix is to prove Proposition \ref{esthom}, Proposition \ref{scsfest} and Lemma \ref{2estspf}
     for $\delta <2/3$. 
     We shall restrict here our analysis to Schr\"odinger Hamiltonians with polynomial potentials,
     for simplicity.  
     We shall use a global pseudo-differential calculus with a  diagonal  metric
     $g = \frac{dx + d\xi}{w(x,\xi)}$ where $w$ is a weight function on the phase space 
     $\R^d_x\times\R^d_\xi$.  The general theory was achieved by L. H\"ormander with the Weyl-Wigner  ordering
      calculus \cite{ho2} after C. Fefferman and R. Beals  for  the "usual" ordering. The construction of parametrices
      for resolvent of elliptic operators is well known. But to cope with $\delta<2/3$  we need to take care of remainder estimates.
      
      \subsection{Parametrix for the resolvent}
      We assume that the potential $V$ is an elliptic polynomial of degree $2k$, which  means that $C_1\<x\>^{2k} \leq V(x) \leq C_2\<x\>^{2k}$
      for $\vert x\vert >R$ where $C>0$, $R>0$.   Denote by 
      $$\hat H(h) = -h^2\Delta + V$$
      and $H(x, \xi) = \vert\xi\vert^2 + V(x)$ its semi-classical symbol. We can get accurate approximations
      for the resolvent $(\hat H(h)-z)^{-1}$ for $h>0$ small and $z\in\C$. It is not difficult to get a formal
      asymptotics as a series in $h$ but to get spectral estimates we need to control the remainder
      terms when~$\vert z\vert$ is large (see~\cite{ro2, daro}).  Recall that the $h$-Weyl quantization of  a  smooth symbol $A$ in $\R^d_x\times\R^d_\xi$ is
       \begin{equation*}
       Op_h(A)\psi(x)= {\hat A}(h)\psi(x) = (2\pi h)^{-d}\iint_{\R^{2d}} A\left(\frac{x+y}{2}, \xi\right)
{\rm e}^{ih^{-1}(x-y)\cdot\xi} \psi(y)dyd\xi.
       \end{equation*}
In particular the Schwartz kernel $K_{A,h}(x,y)$ of ${\hat A}(h)$  is given by
 \begin{equation*}
 K_{A,h}(x,y) = (2\pi h)^{-d}\int_{\R^d}A\left(\frac{x+y}{2}, \xi\right)
{\rm e}^{ih^{-1}(x-y)\cdot\xi}d\xi.
 \end{equation*}
The basic formula for the symbolic computation is the Moyal product formula.
Let $A, B$ be two smooth observables (for example in the Schwartz space ${\cal S}(\R^d\times\R^d)$
and $C$ the $h$-Weyl symbol  of the operator product $\hat C(h):=\hat A(h)\hat B(h)$. Then 
  the $h$-Weyl symbol  of  $\hat C(h)$ is a smooth function
$C(h,x,\xi)$ given by
\beq\label{prod2}
C(h;x,\xi) =  \exp(\frac{ih}{2}\sigma(D_x,D_\xi;D_y,D_\eta))A(x,\xi)B(y,\eta)
\vert_{(x,\xi)=(y,\eta)}, 
\eeq
where $\sigma$ is the standard symplectic form in $ \R^d\times\R^d$, 
$\sigma(x,\xi;xy,\eta) = x\cdot\eta -\xi\cdot y$.\\
    In semiclassical analysis, it is useful to expand 
 the exponent in \eqref{prod2},  so we get the  formal series in $h$:
\bea\label{prod3}
C(h;x,\xi) =  \sum_{j\geq 0}C_j(x,\xi)h^j,\;\; {\rm where}\nonumber\\
C_j(x,\xi) =  \frac{1}{j!}(\frac{i}{2}\sigma(D_x,D_\xi;D_y,D_\eta))^jA(x,\xi)B(y,\eta)
\vert_{(x,\xi)=(y,\eta)}.
\eea
The crucial point here is to have good estimates for the error term not only in $h$ but
also in some extra spectral  parameters as we shall see later
\begin{equation*} 
R_N(h;x,\xi) = C(h;x,\xi) -\sum_{0\leq j\leq N}h^jC_j(x,\xi).
 \end{equation*}
Here one of the symbols $A, B$ is a polynomial in $(x,\xi)$ so the analysis is much simpler.
Let us introduce  the class of symbols $\Sigma(r) $,  $r\in\R$. 
For our application here it is enough to consider the weight function
$$\mu(x,\xi) = (1+\vert x\vert^{2k} + \vert\xi\vert^2)^{1/2k}$$
 (a more general setting is considered in \cite{ro1, daro}).
The condition $A\in\Sigma(r) $ means that for every $j\in\N$ we have
$$
s_j(A,r):= \sup_{\vert\alpha+\beta\vert = j,\; (x,\xi)\in\R^{2d}}
\mu(x,\xi)^{-r +j}\vert \partial_\xi^\alpha\partial_x^\beta  A(x,\xi)\vert < +\infty.
$$
The topology on $\Sigma(r) $ is defined by the semi-norms $s_j(\cdot,r)$. 
A basic result  concerning Weyl quantization is that  for 
every $h>0$, $(A, B)\mapsto A\#B$ is continuous for $\Sigma(r)\times\Sigma(s)$ into $\Sigma(r+s)$.
In  \cite[Appendix~B]{daro},  we can find accurate estimates for $R_N(h;x,\xi)$ in a more general setting.
Here, using that $\hat{H}(h)$  is a polynomial in $(x,\xi)$  we can perform more explicit  computations 
to construct a parametrix for $(\hat H(h)-z)^{-1}$.

Following \cite[p. 134]{ro2}   we construct  a parametrix as follows. By induction on $j\in\N$ we define  
\begin{eqnarray*}
b_{z, 0} &=& (H-z)^{-1} \\
b_{z,j+1} &=& -b_{z, 0} \Big(\sum_{\substack{ \ell +\vert\alpha+\beta\vert = j+1\\  0\leq \ell\leq j  }}\nu(\alpha,\beta)
(\partial_\xi^\alpha\partial_x^\beta H)(\partial_\xi^\beta\partial_x^\alpha b_{z,\ell})\Big)\\
B_{z, N}(h) &=& \sum_{0\leq j\leq N}h^jb_{z,j},
\end{eqnarray*}
where $\nu(\alpha,\beta)=(-1)^{\vert\alpha\vert}[\alpha!\beta!2^{\vert\alpha+\beta\vert}]^{-1}$. Then we have
  \begin{equation}\label{leftp}
  \hat{B}_{z,N}(h)(\hat H(h)-z)=1+\hat{ E}_{z,N}(h),
  \end{equation}
and for our purpose we have to estimate   the error term symbol
$$
E_{z,N}(h) =  B_{z, N}(h)\# ({H}-z) - 1.
  $$        
  This Moyal product  has a finite expansion in $h$ because the symbol  $H$ is a polynomial in $(x,\xi)$.
  Let us recall that    the symbols $b_{z,j}$ have the following properties (see \cite{daro})
   $$
   b_{z, j} = \sum_{0\leq \ell\leq 2j-1}(-1)^{\ell}d_{j\ell}(H-z)^{-\ell-1},\; j\geq 1
   $$
   where $d_{j\ell}$ are polynomials in $\partial_\xi^\alpha\partial_x^\beta H$ for
   $1\leq \vert\alpha+\beta\vert \leq j$. Moreover, for all $m\in \N$ we have 
   if $j=2m$ then $d_{j\ell}=0$  for every $\ell\geq 3m$ and if
   $j=2m+1$ then $d_{j\ell}=0$  for every $ \ell\geq 3m+1$. Furthermore $d_{j\ell}$
   is in the symbol class $\Sigma(2k\ell-2j)$. Thanks to  this vanishing property, we get with the Calderon-Vaillancourt theorem, that   there exists $N_{0}\geq1$   such that for $N\geq 1$
      \begin{equation}\label{normeB}
     \Vert \hat B_{z,N}(h)   \Vert_{{\cal L}(L^2,L^2)}\leq   C_{N} \big(1+h^{N}\vert \Im z\vert^{-3N/2}\big)\vert\Im z\vert^{-N_0}. 
   \end{equation}
Next, with  \eqref{prod3} we get the  following estimate:
   There exists $n_0$ such that for every $N\geq 1$, $(\alpha, \beta)\in\N^{2d}$,  there exists $C>0$ such that 
   for every $h\in]0, 1]$ and $z\in\C_{\gamma_0}$, we have 
   \bea\label{acres1}
   \vert\partial_\xi^\alpha\partial_x^\beta  E_{z,N}(h) \vert \leq
   Ch^{N+1}\sum_{\frac{N+1}{2k}\leq \ell\leq 3N/2+n_0+\vert\alpha+\beta\vert}
   \left\vert\frac{H}{H-z}\right\vert^{\ell+1}
     ( \mu(x,\xi))^{-2(N+1)-\vert(\alpha+\beta)\vert}.
  \eea
  In the r.h.s of \eqref{acres1}  it is desirable to have  large  decay in $\vert z\vert$ and in $\mu$. First, it  is known that there exists $\gamma_0\in \R$ such that the spectrum of $\hat H(h)$ is in $[\gamma_0, +\infty[$ for all $h\in ]0, 1]$.    Then, with elementary considerations, we have that  for every $\theta\in[0, 1]$   
  $$
   \left\vert\frac{H}{H-z}\right\vert \leq \frac{1}{\vert\sin\phi\vert} \vert z\vert^{-\theta} H^{\theta}\leq  \frac{1}{\vert\sin\phi\vert} \vert z\vert^{-\theta} \mu^{2k\theta},
    $$
   where $\phi =\arg z$, $0 <\vert\phi\vert\leq \pi/2$, $\vert z\vert \geq \gamma_0/2$.\\
   So we can choose $\theta\in]0, 1[$  such that there exist  positive numbers 
   $\delta_1, \delta_2, \delta_3,\delta_{4}>0$, $N_{0}>0$  and $C>0$ such that for $\vert z\vert \geq \gamma_0/2$,
   $0<\vert \phi\vert\leq\pi/2$, $(x,\xi)\in\R^{2d}$, $N\geq 1$,  we have
   \beq\label{param1}
    \vert\partial_\xi^\alpha\partial_x^\beta  E_{z,N}(h) \vert \leq
   Ch^{N+1}\vert\sin\phi\vert^{-3 N/2-N_{0}-\vert\alpha+\beta\vert}\vert z\vert^{-\delta_1N-\delta_2\vert\alpha+\beta\vert}
   \mu^{-\delta_3N-\delta_4\vert\alpha+\beta\vert}.
   \eeq
   As a result, by the Calderon-Vaillancourt theorem, we obtain that  $\hat E_{z,N}(h)$ is continuous from $L^2(\R^d)$ into ${\cal H}_k^{\delta N}$  for some $\delta >0$ and that there exist $N_{0}\geq1$, $\kappa>0$, $C>0$  such that for $N\geq 1$
   \begin{equation}\label{normeE}
     \Vert \hat E_{z,N}(h)   \Vert_{{\cal L}(L^2,{\cal H}_k^{\delta N})}\leq   Ch^{N+1}\vert\sin\phi\vert^{-3 N/2-N_{0}}\vert z\vert^{-\kappa N}. 
   \end{equation}
     Since we work with the Weyl quantization, we can get a  right parametrix by taking the adjoint of\,\eqref{leftp}, and therefore 
       \begin{equation*}
(\hat H(h)-z)  \hat{B}^*_{\ov{z},N}(h)=1+\hat{ E}^*_{\ov{z},N}(h),
  \end{equation*}
    which combined with \eqref{leftp} yields
  \beq\label{param2}
  \hat B_{z,N}(h) - (\hat H(h)-z)^{-1} = \hat E_{ z,N}(h)\hat B^*_{\bar z,N}(h) -
  \hat E_{z,N}(h)(\hat H(h)-z)^{-1} \hat E^*_{\bar z,N}(h) := \hat R_{z,N}(h).
  \eeq
   So we can  easily see that $\hat R_{z,N}(h)$ is a smoothing operator in the weighted Sobolev
   scale spaces $\{{\cal H}_k^s\}_{s\in\R}$ when $N$ is large, and we can use the next lemma  proved in \cite[Proposition 1.3]{ro1} to  give an estimate of the kernel of $\hat R_{z,N}(h)$.
   \begin{lemm}\ph \label{estProd} Let  $\hat R, \hat S$ be bounded operators in $L^2(\R^d)$ and assume that 
   $\hat R(L^2(\R^d))\subseteq {\cal H}_k^m$ and  $\hat S(L^2(\R^d))\subseteq {\cal H}_k^m$
   with $m>\frac{kd}{k+1}$. Then the operator $\hat T := \hat R\hat S^*$ has a  (Lipschitz) continuous  Schwartz kernel $K_T(x,y)$
   on $\R^d\times\R^d$ and we have the estimate 
   \beq\label{param3}
   \vert K_T(x,y)\vert \leq C_{k,d,m}W_m(x)W_m(y)\Vert\hat R\Vert_{{\cal L}(L^2,{\cal H}_k^m)}
   \Vert\hat S\Vert_{{\cal L}(L^2,{\cal H}_k^m)}
   \eeq
   where $C_{k,d,m}$ only depends  on $k,d,m$ and $W_m(x) = (1+\vert x\vert^{2k})^{d/2-(k+1)m/(2k)}$.
   \end{lemm}
 \begin{rema}
 If $m$ can be taken larger then we get that $K_T$ has derivatives in $(x,y)$ with corresponding estimates.
 \end{rema}

Using \eqref{param1}, \eqref{param2}, \eqref{param3}, we get  the following estimate for the error term in the parametrix. Let us denote by $K_{R_{z,N}(h)}$ the Schwartz kernel of $\hat R_{z,N}(h)$.
\begin{lemm}\ph
There exist $N_0$ large enough and $\kappa_{1},\kappa_{2}>0$, $C>0$  such that for $N\geq 1$
\beq\label{erparam}
K_{R_{z,N}(h)}(x,y) \leq Ch^{N+1}\big(1+h^N \vert  \sin \phi \vert^{-3 N/2}\big) \vert  \sin \phi \vert^{-3 N/2-N_0}|z|^{-\kappa_{1}N}(1+\vert x\vert^{2k}+\vert y\vert^{2k})^{-\kappa_{2}N}
\eeq
for all $h\in]0, 1]$, $|z|\geq \gamma_{0}/2$ with $ \arg z =\phi$  and $x, y\in\R^d$.

Moreover, if $z\in \C\backslash\R$ is such that $\vert z \vert\leq A$, 
\beq\label{erparam2}
K_{R_{z,N}(h)}(x,y) \leq Ch^{N+1}\big(1+h^N \vert  \Im z \vert^{-3 N/2}\big) \vert  \Im z \vert^{-3 N/2-N_0}(1+\vert x\vert^{2k}+\vert y\vert^{2k})^{-\kappa_2N}.
\eeq
\end{lemm}

 \begin{proof} The estimate of \eqref{erparam2} is a direct consequence of \eqref{erparam} by writing $|\sin \phi|\geq |\Im z|/A$. For\,\eqref{erparam}, we use Lemma \ref{estProd} twice, using the estimates \eqref{normeB} and \eqref{normeE}. We write  $\hat E_{ z,N}(h)\hat B^*_{\bar z,N}(h) =\hat R \hat S^*$ with $\hat R=\hat E_{ z,N}(h) (\hat H(h)-z)^{\eta N} $ and ${\hat S=\hat B_{\ov  z,N}(h) (\hat H(h)-\ov z)^{-\eta N} }$, with $\eta>0$ small enough so that the hypothesis of Lemma \ref{estProd} is satisfied. For the second term, we write  \\
${E_{z,N}(h)(\hat H(h)-z)^{-1} \hat E^*_{\bar z,N}(h)=\hat R \hat S^*}$ with $\hat R=\hat E_{ z,N}(h) $ and ${\hat S=\hat E_{\ov  z,N}(h) (\hat H(h)-\ov z)^{-1} }$.
  \end{proof}
 
\subsection{Proof of Proposition \ref{esthom}}\label{appA2}
  \begin{proof}
   We   use the previous result with $h=1$, and we denote by $\hat H=\hat H(1)$.  The heat operator ${\rm e}^{-t\hat H}$ is related with the resolvent by a Cauchy integral
   \begin{equation*}
      {\rm e}^{-t\hat H} = \frac{1}{2i\pi}\oint_{\Gamma}{\rm e}^{-tz}(\hat H -z)^{-1}dz
   \end{equation*}
   where $\Gamma$ is the following contour in the complex plane. Fix $\phi\in]0, \pi/2]$. 
   Let $z_0=(\gamma_0/2){\rm e}^{i\phi}$. Let $\Gamma_+$  be the line $[z_0+r{\rm e}^{i\phi}, r\geq 0]$,
   $ \Gamma_- =\overline{\Gamma_+}$, $\Gamma_0 =[(\gamma_0/2){\rm e}^{i\psi}, -\phi\leq\psi\leq \phi]$.
   So $\Gamma =\Gamma_0\cup\Gamma_-\cup\Gamma_+$ with a suitable orientation.
   
   We give the main steps of the end of the proof. 
   \begin{itemize}
   \item[$\bullet$]  Thanks to the parametrix computed for $(\hat H -z)^{-1}$, we have
   \begin{equation*}
      {\rm e}^{-t\hat H} = \frac{1}{2i\pi}Op_1\Big(\sum_{j=1}^{2N}f_{j}(x,\xi)\oint_{\Gamma}{\rm e}^{-tz}(H(x,\xi) -z)^{-j}dz\Big)- \frac{1}{2i\pi}  \oint_{\Gamma}{\rm e}^{-tz}\hat {R}_{z,N}(1)dz ,
   \end{equation*}
   where the $f_{j}$ are linear combinations of the $d_{j\ell}$, hence they are polynomials in $(x,\xi)$.
      \item[$\bullet$] The kernel of the main contribution is obtained with the 
residue theorem: 
   \begin{equation*}
    \frac{1}{2i\pi}\oint_{\Gamma}{\rm e}^{-tz}(H(x,\xi) -z)^{-j}dz=\frac{t^{j-1}}{(j-1) !}\e^{-tH(x,\xi)},
   \end{equation*}
   and using degree considerations of the $f_{j}$.
        \item[$\bullet$] For $N\geq1$ large enough, the kernel of the remainder term is estimated with \eqref{erparam}.
        \end{itemize}
   \end{proof}
   
   
   \subsection{Proof of Proposition \ref{scsfest} }\label{A.3}
   \begin{proof}
    We use here the semi-classical  functional calculus for smooth functions with compact support, and we apply it to the operator $\hat H(h)=-\h^{2} \Delta +V_{\eps}$ where $V_{\eps}(x)=\eps^{2k} V(x/\eps)$.  We can check that all the previous estimates hold true uniformly in $\eps>0$.
    
     Let $f$ be a  non negative $C^\infty$ function in $]-2C_0, 2C_0[$ with  a compact support,
    such that $f=1$ in $[-C_0,  C_0]$.  Using the spectral theorem for general self-adjoint operators,
    we can consider the new operator 
   \beq\label{sm}
    g\big(\hat H(h)\big) =f\left(\frac{\hat H(h)-\nu}{h^{\delta}}\right).
  \eeq
  The operator  $g\big(\hat H(h)\big) $ has clearly a smooth Schwartz kernel  $K_{f,h}$ (compute it in a basis of  eigenvector of~$\hat H(h)$)
   and we  have
\begin{equation*}   
   \vert\pi_{H(h)}(\nu+\mu h;x,x) - \pi_{H(h)}(\nu;x,x)\vert \leq  K_{f,h}(x,x),\;\; \forall x\in\R^d,
    \end{equation*}
    where $|\mu|\leq C h^{\delta-1}$.  So it is enough to show that $K_{f,h}(x,x)={\cal O}(h^{\delta-d})$ uniformly in $x\in\R^d$.
     
     Let us recall the almost-analytic formula for the functional calculus (see \cite[Chapter 8]{disj}) 
     \beq\label{anfu}
     g(\hat H(h)) = \frac{1}{\pi}\int_\C\bar\partial\tilde g(z)(\hat H(h)-z)^{-1}L(dz),
\eeq 
where $L(dz)=dxdy$ is the Lebesgue measure and     where $\tilde g$ is an almost analytic extension of $g$.
     This means that $g$ is $C^\infty$ in $\C$ and satisfies, for all $M\geq 0$, 
\begin{equation*}     
     \vert\bar\partial\tilde g(z)\vert \leq  C_M\vert\Im z\vert^M,\;\; \forall z\in \C.
    \end{equation*}
      We can choose $\tilde g$ supported 
      in a small complex neighborhood of $]-2C_0, 2C_0[$.  Moreover it results from~\cite{disj}
     that $C_M$  is estimated by some semi-norms of $g$. 
      \beq\label{aana2}
     \vert\bar\partial\tilde g(z)\vert \leq
       \gamma _M\vert\Im z\vert^M\big(\Vert\widehat{D_x^{M+1} g}\Vert_{L^1} + \Vert D_x^{M+2} g\Vert_{L^1} \big),\;\;\; \forall z\in \C.
    \eeq
    
      We give the main steps of the end of the proof. 
   \begin{itemize}
   \item[$\bullet$]  We plug the parametrix of $(\hat H(h)-z)^{-1}$
 in \eqref{anfu}  to compute the kernel of  $g(\hat H(h))$.
       \item[$\bullet$] The kernel of the main contribution $\hat B_{z,N}$ is computed explicitly. Actually, the first term in the expansion of  $g(\hat H(h))$ is 
       \begin{equation*}
       Op_{h}\Big(\frac{1}{\pi}\int_\C\bar\partial\tilde g(z)(  H-z)^{-1}L(dz)\Big)= Op_{h}\big(g(H) \big)=    Op_{h}\Big(   f\Big(\frac{  H-\nu}{h^{\delta}}\Big)  \Big),
       \end{equation*}
where we used  that $\frac{1}{\pi z}$  is a fundamental solution of~$\bar\partial$. Therefore, the principal term of $K_{f,h}$ is given by the following formula
  \beq\label{pterm}
   K^0_{f,h}(x,x) = (2\pi h)^{-d}\int_{\R^d} f\left(\frac{H(x,\xi)-\nu}{h^{\delta}}\right)d\xi,
   \eeq
which implies $K^{0}_{f,h}(x,x)={\cal O}(h^{\delta-d})$.
          \item[$\bullet$] Choosing $M$ large enough we can estimate the contribution
   of $\hat R_{z,N}(h)$ by using \eqref{erparam2} and~\eqref{aana2} for the smallest  integer $M \geq 3N+N_0$.
 In  \eqref{aana2} the loss  in $h$ is $h^{-\delta M}$  but this is compensated
   by the factor~$h^{N+1}$ for $N$ large if we choose $\delta < 2/3$. 
        \end{itemize}
     \end{proof}

  \subsection{Proof of  Lemma \ref{2estspf}} \label{A.4}
  \begin{proof}
     The proof uses the same tools as for the proof of Proposition \ref{scsfest}.
     We choose two cutoff functions $f_\pm$ with $f_+$ as above and $f_-$ such
     that supp$(f_-)\subseteq]C_1, C_0[$, $f_-=1$ in $[2C_1, C_0/2]$ where $C_1 < C_0/4$.
     If $K_{f,h}$ is the Schwartz kernel  of $g(\hat H(h))$ given par \eqref{sm} we have
\begin{equation*}
     K_{f_-,h} \leq e_{x,h} \leq K_{f_+,h}(x).
\end{equation*}
     So we have to prove
     \beq\label{reduc}
  C_0N_hh^{\beta_{2p,\theta}}\leq \left(\int_{\R^d}\<x\>^{k\theta(p-1)}K^{p}_{f,h}(x)dx\right)^{1/p}
  \leq   C_1N_hh^{\beta_{2p,\theta}}, 
  \eeq
   for  every $f$ like in \eqref{sm}. 
   
     We use again the almost-analytic functional  calculus. From \eqref{erparam} we see that
   the contribution of the error term of the resolvent parametrix is  of order ${\cal O}(h^{N+1})$. As previously, the principal term of~$K_{f,h}$ is given by \eqref{pterm}, hence we get \eqref{reduc}    from a two  side estimate of this integral.
     \end{proof}  
 
\subsection{On the case $\delta=1$}
Here we show how Proposition \ref{scsfest}  for $\delta =1$  can be deduced   directly from results established in \cite{KTZ}. We  admit here that these results can be extended for $V$ depending on a parameter $\eps$.\\[2pt]
\indent $\bullet$ {\bf Estimates outside the turning points:} Outside the turning points $V(x) = \nu$ the estimate can be proved using a standard WKB approximation
 for the propagator $U(t):={\rm e}^{-ith^{-1}\hat H(h)}$. 
 Let   $U(t,x,y)$ be the Schwartz kernel of $U(t)$.  Let us   give here a sketch of the proof.\\
 Let $\rho\in C_0^\infty   \big(\,]-T_0, T_0[\,\big)$  with $T_0>0$ small enough. We have
 $$
 I_{x,\nu}(h):= \int\rho(t){\rm e}^{it\nu h^{-1}}U(t,x,x)dt = \sum_{j\geq 0}\hat\rho\left(\frac{\lambda_j-\nu}{h}\right)\vert\varphi_j(x)\vert^2.
 $$
 Choosing $\rho$ even such that $\hat\rho\geq 0$ and $\hat\rho(0)=1$ it is enough to prove
 \beq\label{NT}
  I_{x,\nu}(h) = {\cal O}(h^{1-d}).
  \eeq
  We consider now a WKB approximation for   $I_{x,\nu}(h)$
  $$
  I_{x,\nu}(h) \approx (2\pi h)^{-d}\int_\R \int_{\R^d}\rho(t){\rm e}^{ih^{-1}(S(t,x,\eta)-y\cdot \eta +t\nu)}
  \left(\sum_{j\geq 0} h^ja_j(t,x,\eta)\right)dtd\eta.
  $$
  
  Using a localization energy argument it is enough to consider a bounded  open set of the phase space
 ${\cal V} =H^{-1}]\nu_0-\eps_0, \nu_0 +\eps_0[$.  \\
 The term  $S(t,x,\eta)$ is the solution of the Hamilton-Jacobi equation in ${\cal V}$ and  the $a_j$
 are determined by transport equations (see \cite{ro2} for details).\\
 The stationary points of the phase $\Phi_{x,\nu}(t,\eta) : = S(t,x,\eta)-y\cdot \eta +t\nu$ 
 satisfy the equations $t=0$, $\vert\eta\vert^2 = \nu -V(x)$ if $T_0$ is chosen small enough.
 So  if $V(x)<\nu$ the critical set is a smooth submanifold~${\cal C}_{x,\nu}$ of $\R\times\R^d$
 of codimension 2 and the Hessian of $\Phi_{x,\nu}$ is non degenerate in the normal directions
 to~${\cal C}_{x,\nu}$.  So the stationary phase theorem gives
 that $I_{x,\nu}(h) = {\cal O}(h^{1-d})$.\\
 If $V(x) >\nu$ the non-stationary phase theorem gives that $I_{x,\nu}(h) = {\cal O}(h^{\infty})$.\\[2pt]
 \indent $\bullet$ {\bf Estimates at the turning points:} Now we consider the case $\vert V(x) -\nu\vert <\eps_1$ for $\eps_1>0$ arbitrary small. We give the argument  used in \cite{KTZ}, which  is completely different.\\
 Let $\Pi_h$ the spectral projector for $\hat H(h)$ on $[\nu, \nu+\mu h]$. We claim that it is enough to prove
 \beq\label{Proj}
 \Vert\Pi_h\Vert_{L^2\to L^\infty} = {\cal O}(h^{\frac{1-d}{2}}).
 \eeq
 Assume that \eqref{Proj} holds true. Denote by  $e_h(x,y)$  the Schwartz kernel of the projector $\Pi_h$ and set $u_{x}: y \mapsto e_{h}(x,y)$. Then  we have $\Pi_{h} u_{x}=u_{x}$, thus  
 $$
e_{h}(x,x)= u_{x}(x) = (\Pi_h u_{x})(x) \leq \Vert\Pi_h\Vert_{L^2\to L^\infty}\Vert u_x \Vert_{L^2}.
 $$ 
Then using that $\Pi^{2}_{h}=\Pi_{h}$ we get 
\begin{equation*}
\Vert u_x \Vert_{L^2}^{2}=\int \vert e_{h}(x,y)\vert^{2}d y =e_{h}(x,x),
\end{equation*}
which entails that $e_h(x,x)={\cal O}(h^{1-d})$, thanks to \eqref{Proj}.

 Let $\chi$ be supported in ${\cal V}\cap\{\vert V(x) -\nu\vert <\eps_1\}$. 
 It is enough to prove that 
 \beq\label{KTZ}
 \Vert\hat\chi(h)\Pi_h\Vert_{L^2\to L^\infty} = {\cal O}(h^{\frac{1-d}{2}}).
 \eeq
Using that 
$\Vert(\hat H(h)-\nu)\Pi_h\Vert_{L^2\to L^2} = {\cal O}(h)$, for $d\geq 3$ estimate \eqref{KTZ} is a  direct consequence
 of  \cite[Theorem 6]{KTZ}. For $d=2$  this theorem gives the estimate 
 ${\cal O}(h^{-1/2}\vert\log\h\vert^{1/2})$, but applying the more difficult   \cite[Theorem 6]{KTZ}, the $\log$ term
 can be eliminated.

\end{document}